\documentclass[11pt]{article}
\usepackage{amssymb}
\usepackage{amsmath}
\usepackage{amsthm}
\usepackage{latexsym}
\usepackage{amsfonts}
\usepackage{graphicx}
\usepackage{graphics}

\newcommand{\dis}{\displaystyle}
\theoremstyle{plain}
\newtheorem{thm}{Theorem}[section]   
\newtheorem{prop}{Proposition}[section] 
\newtheorem{cor}{Corollary}[section] 
\newtheorem{Def}{Definition}[section] 
\newtheorem{lem}{Lemma}[section] 

\theoremstyle{definition}
\newtheorem{rem}{Remark}[section]

\newtheorem{obs}{Observation}[section] 

\newtheorem*{Proof}{Proof}

\newcommand{\fa}{\forall}

\newcommand{\el}{\ell}

\newcommand{\ra}{\;\rightarrow\;}

\newcommand{\al}{\alpha}
\newcommand{\bi}{\beta}
\newcommand{\ga}{\gamma }
\newcommand{\Ga} {{\varGamma}}

\newcommand{\de}{\delta }

\newcommand{\De} {{\varDelta}}
\newcommand{\e}{\varepsilon }
\newcommand{\f}{\varphi}
\newcommand{\Fi}{\varPhi}

\newcommand{\thi}{\theta }
\newcommand{\vthi}{\vartheta }
\newcommand{\Thi} {{\mit\Theta}}
\newcommand{\X} {{\mit\Xi}}
\newcommand{\La} {{\mit\Lambda}}

\newcommand{\la}{\lambda }
\newcommand{\mi}{\mu }

\newcommand{\si}{\sigma }
\newcommand{\ups}{\upsilon}
\newcommand{\ti}{\tau }

\newcommand{\R}{\mathbb{R}}

\newcommand{\N}{\mathbb{N}}

\newcommand{\ssum}{\sum\limits}

\newcommand{\ld}{\ldots}

\newcommand{\qb}{$\quad\blacksquare$}
\newcommand{\mf}{\mathfrak{M}}

\begin{document}
\title{\bf Simultaneous generic approximation by the iterates of the Cesaro operator}
\author{Th. Douvropoulos}
\date{}
\maketitle
\begin{abstract}
We show that for the generic sequence $a$ of elements in a subset
$A$ of a separable locally convex metrizable space $V$, the
sequences $[T^k(a)]_n$, $n=1,2,\ld$ are dense in the convex hull
conv$A$ of $A$ for all $k=1,2,\ld$; where $T$ is the Cesaro
operator. Further, if conv$A$ is dense in $V$, then for every
sequence $x_k\in V$, $k=1,2,\ld$ there exists $\la_n\in\{1,2,\ld\}$,
$n=1,2,\ld$, such that, we have the simultaneous approximation
$[T^k(a)]_{\la_n}\!\!\ra\!\! x_k$, as $n\!\!\ra\!\!+\infty$ for all
$k=1,2,\ld\;.$ These phenome-\linebreak na are topologically generic
and in the case where $A$ is a vector space they are algebraically
generic.
\end{abstract}
{\em 2010 AMS classification number}: 41A36, 41A28, 40A05. \vspace*{0.1cm} \footnotetext{This work was completed as part of
the implementation of a Master's program that was partly funded by
the Act: ``SSF scholarships granting Program, through the process of
individualized evaluation for the academic year 2011-2012'' by
resources of the Lifelong Learning Programme (LLP) of the European
Social Fund (ESF) and of the NSRF 2007-2013.}
\section{Introduction}\label{sec1}
\noindent

In 1896 Emile Borel in a letter to Mittag Leffler published in Acta
Mathematica claimed that in ``general'' the circle of convergence of
a power series is a natural boundary. In \cite{4} Baires Category
theorem is used to prove that this holds generically for every
holomorphic function in the open unit disc. This Baire's method has
been used extensively in analysis; see \cite{3}, \cite{2}, as well
\cite{5}, \cite{6} and \cite{1}. During the Congress of the Hellenic
Mathematical Society (November 2011) V. Nestoridis illustrated the
Baire's method in a very simple example, related to the Cesaro
operator.

If $a=(a_n)^\infty_{n=1}$ is a sequence of real numbers then
$[T(a)]_{n}=\frac{a_1+\cdots+a_n}{n}$. It is well known that if
$a_n\ra\el$, $\el\in\R$ as $n\ra+\infty$ then $[T(a)]_{n}\ra\el$ as
well, as $n\ra+\infty$. The converse is not true; for instance if
$a_n=(-1)^n$, then $[T(a)]_{n}\ra0$, as $n\ra+\infty$ but $a_n$ is
not convergent. Thus, if we start by an arbitrary sequence $a$
(which in general does not converge) we wish to apply an iterate of
the Cesaro operator $T:\R^\N\ra\R^\N$ say $T^k$ defined on the space
$\R^\N$ of all sequences of real numbers and we hope that for some
$k\in \N$ the sequence $[T^k(a)]_{n}$, $n=1,2,\ld$ becomes a
convergent sequence. Unfortunately this is not always possible. More
precisely for a large subset $X$ of $\R^\N$ (a $G_\de$ and dense
subset) the sequence $[T^k(a)]_{n}$, $n=1,2,\ld$ is dense in $\R$
and this simultaneously for all $k=1,2,\ld\;.$

The previous set $X$ (possibly empty) has the following description
\[
X=\bigcap_{k,j,s}\;\;\bigcup^\infty_{n=1}\bigg\{a\in\R^N:|[T^k(a)]_{n}-q_j|<\frac{1}{s}\bigg\}
\]
where $q_j$ is an enumeration of $Q$, $s=1,2,\ld$ and $k=1,2,\ld\;.$

Because the operator $T:\R^\N\ra\R^\N$ is continuous, where $\R^\N$
is endowed with the cartesian topology (and therefore $\R^\N$ is a
complete metric space) it follows that the sets
$\Big\{a\in\R^\N:|[T^k(a)]_{n}-q_j|<\frac{1}{s}\Big\}$ are open for
all $k,j,s$ and $n$ and therefore the set $X$ is $G_\de$. In order
to apply Baire's Category theorem and prove that $X$ is $G_\de$ and
dense and therefore $X\neq\emptyset$, it suffices to prove that
$\bigcup\limits^\infty_{n=1}\Big\{a\in\R^\N:|[T^k(a)]_{n}-q_j|<\frac{1}{s}\Big\}$
is dense for all $k$, $j$ and $s$. Towards this end we consider
$M\in\N$ and $\ga_1,\ld,\ga_M$ arbitrary real numbers and we are
looking for a sequence $a\in \R^N$ and a natural number $n\in\N$ so
that $a_i=\ga_i$ for $i=1,\ld,M$ and $[T^k(a)]_{n}-q_j|<\frac{1}{s}$
(we remind that $k,j,s$ are fixed). We consider the sequence $a$
given by $a_i=\ga_i$ for $i=1,\ld,M$ and $a_i=q_j$ for all $i\ge
M+1$. This sequence converges to $q_j$; therefore
$\big[T^k(a)\big]^\infty_{n=1}$ converges also to $q_j$. It follows
that there exists $n\in\N$ such that
$[T^k(a)]_{n}-q_j|<\frac{1}{s}$, and the proof is complete.

After this simple argument V. Nestoridis asked for an explicit
example of such a sequence $a\in X$. Towards this let $q_j$ be a
standard enumeration of the set $Q$; then $|q_j|\le j$. We consider
the sequence $a$ given by
\[
q_1,\underbrace{q_2,\ld,q_2}_{k_2\;\text{times}},\ld,\underbrace{q_n,\ld,q_n}_{k_n\;\text{times}}
,\ld,
\]
where $k_n=n^{n^3}$. Then one can easily check that $a\in X$.

Further, if we start with any subset $A\subset\R$ (or of any
separable locally convex metrizable vector space) we show that there
are sequences $a$ with elements in $A$, such that $[T^k(a)]_n$,
$n=1,2,\ld$ is dense in the convex hull conv$A$ of $A$ and this for
all $k=1,2,\ld\;.$ The set $X$ of such sequences is $G_\de$ and
dense in $A^N$. This is established in \S\,3 of the present paper
and the proof uses Baire's Category theorem. This is done for the
sake of simplicity, while an other proof without using Category is
also possible. We also mention that the previous set $X$ contains a
dense vector subspace in $A^N$ except 0, provided that $A$ is a
vector space. This is not proved in \S\,3, but it follows from the
results of \S\,4 (see Remark \ref{rem4.11}).

Finally, since for $a\in X$ the sequences $[T(a)]_n$, $[T^2(a)]_n$,
$n=1,2,\ld$ are dense in conv$A$ it follows that for any
$x_1,x_2\in$conv$A$ there exist two sequences of natural numbers
$\al_n$, $\mi_n$, $n=1,2,\ld$ so that $[T(a)]_{\la_n}\ra x_1$ and
$[T^2(a)]_{\mi_n}\ra x_2$, as $n\ra+\infty$. The question that
naturally arises is whether we can have $\la_n=\mi_n$; that is, is
it possible to have simultaneous approximation?

We answer in the affirmative the previous question in \S\,4 provided
that conv$A$ is dense. The result is topologically generic and in
the case where $A$ is a dense vector space, then the result is
algebraically generic. However we mention that in general the
simultaneous approximation by the iterates of the Cesaro operator is
not possible; see Proposition \ref{prop4.12}. We also mention that
in the proof of \S\,4 we avoid to use Baire's theorem.

It would be interesting to know, if for other natural operators such results hold, or if the family $\{T^k:k\ge1\}$ could be replaced by other families of operators.
\section{Preliminaries}\label{sec2}
\noindent

For an arbitrary vector space $V$, we consider the Cesaro operator
$T:V^\N\ra V^\N$ such that if $\vthi\equiv(\vthi_n)\in V^\N$, then
$[T(\vthi)]_n=\frac{\vthi_1+\cdots+\vthi_n}{n}$. For each
$k=2,3,\ld$ we define $T^k=T\circ T^{k-1}$.

Next, we consider  the matrix that describes $T^k$. For
$n,m\in\N^\ast$, let $T^k_{(n,m)}$ be its element of the
$n^{\text{th}}$ row and $m^{\text{th}}$ column.

For instance the matrix that describes $T$ is:
\[
T=\left(\begin{array}{ccccccc}
          1 & 0 & 0 & 0 & \cdot & \cdot & \cdot \\
          \frac{1}{2} & \frac{1}{2}& 0 & 0 & \cdot & \cdot & \cdot \\
          \frac{1}{3} & \frac{1}{3} & \frac{1}{3} & 0 & \cdot & \cdot & \cdot \\
          \cdot & \cdot & \cdot & \cdot & \cdot & \cdot & \cdot \\
          \cdot & \cdot & \cdot & \cdot & \cdot & \cdot & \cdot \\
          \frac{1}{n} & \frac{1}{n} & \cdots & \frac{1}{n} & \cdot & \cdot & \cdot \\
          \cdot & \cdot & \cdot & \cdot & \cdot & \cdot & \cdot \\
          \cdot & \cdot & \cdot & \cdot & \cdot & \cdot & \cdot
        \end{array}\right)
\]
For the elements $T^k_{(n,k)}$ of the matrix of $T^k$
$(k\in\N^\ast)$ we prove the following:
\setcounter{prop}{-1}
\begin{prop}\label{prop2.0}
The matrix of $T$ is lower-triangular and the same is true for the
matrix of $T^k$ since it is a product of lower-triangular matrices.

Moreover the sum of the elements of the $n^{\text{th}}$ row of the
matrix of $T$ is equal to 1. This is also true for $T^k$
$k\in\N^\ast)$.
\end{prop}
\begin{Proof}
Indeed, since it is by definition true for $k=1$, we suppose it is
true for $k\in\N^\ast$ and we only need to prove it is true for
$k+1$.

Let $n\in\N^\ast$:
$\ssum^\infty_{m=1}T^{k+1}_{(n,m)}=\ssum^\infty_{m=1}\Big(\ssum^\infty_{i=1}T_{(n,i)}\cdot
T^k_{(i,m)}\Big)=\ssum^\infty_{i=1}\Big(\ssum^\infty_{m=1}T_{(n,i)}\cdot
T^k_{(i,m)}\Big)=\ssum^\infty_{i=1}\Big(T_{(n,i)}\cdot\Big(\ssum^\infty_{m=1}T^k_{(i,m)}\Big)\Big)$
which by use of the inductive argument becomes:\\
$=\ssum^\infty_{i=1}T_{(n,i)}=1$. \qb
\end{Proof}
\begin{prop}\label{prop2.1}
Let $k,n\in\N^\ast$ be given. Then $T^k_{(n,m)}$ is a decreasing sequence with respect to $m\in\N^\ast$.
\end{prop}
\begin{Proof}
For $k=1$, $T_{(n,m)}=\frac{1}{n}$ when $m\le n$ and $T_{(n,m)}=0$
when $m>n$. This is indeed decreasing with respect to $m$. Suppose
that the proposition is true for $k\in\N^\ast$. We will show that it
is true for $k+1$:

Let $k,n,m\in\N^\ast$:
\begin{align*}
T^{k+1}_{(n,m)}-T^{k+1}_{(n,m+1)}&=\sum^\infty_{i=1}(T_{(n,i)}\cdot T^k_{(i,m)}-T_{(n,i)}\cdot T^k_{(i,m+1)}) \\
&=\sum^\infty_{i=1}[T_{(n,i)}\cdot (T^k_{(i,m)}-T^k_{(i,m+1)})]\ge0,
\end{align*}
since by the inductive argument $T^k_{(i,m)}\ge T^k_{(i,m+1)}$
$\fa\; i=1,2,\ld\;.$  \qb
\end{Proof}
\setcounter{cor}{1}
\begin{cor}\label{cor2.2}
Let $k,n\in\N^\ast$, then if $m>\frac{n}{2}$, $m\in\N$ it is true
that $T^k_{(n,m)}\le\frac{2}{n}$.
\end{cor}
\begin{Proof}
By Proposition \ref{prop2.0} $\ssum^\infty_{m=1}T^k_{(n,m)}=1$ and
by Proposition \ref{prop2.1} $T^k_{(n,m)}$ is decreasing with
respect to $m$. Thus if for some $m>\frac{n}{2}$,
$T^k_{(n,m)}>\frac{2}{n}$ then $T^k_{(n,i)}>\frac{2}{n}$ for each
$i=1,2,\ld,\lceil\frac{n}{2}\rceil$, thus
$\ssum^\infty_{i=1}T^k_{(n,i)}>\frac{n}{2}\cdot\frac{2}{n}=1$ which
is a contradiction.  \qb
\end{Proof}
\setcounter{prop}{2}
\begin{prop}\label{prop2.3}
For each $k,n,m\in\N^\ast$, it is true that:
\[
T^k_{(n,m)}\le\frac{1}{n}\bigg[1+\frac{1}{1!}\cdot(\log n)^1+\cdots+\frac{1}{(k-1)!}(\log n)^{k-1}\bigg].
\]
\end{prop}
\begin{Proof}
For $k=1$, the above inequality, becomes $T_{(n,m)}\le\frac{1}{n}$ which is true since $T_{(n,m)}=\frac{1}{n}$ or 0.

We suppose the statement is true for $k\in\N^\ast$. For $k+1$: $T^{k+1}_{(n,m)}=\ssum^\infty_{i=1}T_{(n,i)}\cdot T^k_{(i,m)}$. But the matrices of $T$ and $T^k$ are lower-triangular thus $=\ssum^n_{i=m}T_{(n,i)}\cdot T^k_{(i,m)}=\frac{1}{n}\cdot\ssum^n_{i=m}T^k_{(i,m)}$ which by the inductive argument is
\begin{eqnarray}
\le\frac{1}{n}\cdot\sum^n_{i=m}\frac{1}{i}\bigg[1+\frac{1}{1!}\cdot(\log i)^1+\cdots+\frac{1}{(k-i)!}\cdot(\log i)^{k-1}\bigg].  \label{sec2eq1}
\end{eqnarray}
However if $x\in\R$, $x\ge1$ then $f(x)=\frac{1+\log
x+\cdots+\frac{1}{(k-1)!}\cdot(\log x)^{k-1}}{x}$ is a decreasing
function.

Indeed
\begin{align*}
f'(x)=&\;\frac{\Big[\frac{1}{x}+\frac{1}{1!}(\log x)^2+\cdots+\frac{1}{(k-2)!}(\log x)^{k-2}\cdot\frac{1}{x}\Big]\cdot x}{x^2} \\
&-\frac{\Big[1+\frac{1}{1!}(\log x)^2+\cdots+\frac{1}{(k-1)!}(\log x)^{k-1}\Big]}
{x^2}\\
=&\;\frac{\frac{-1}{(k-1)!}(\log x)^{k-1}}{x^2}<0.
\end{align*}
Thus
\begin{align*}
&\sum^n_{i=m}\frac{1}{i}\bigg[1+\frac{1}{1!}(\log i)^1+\cdots+\frac{1}{(k-1)!}(\log i)^{k-1}\bigg] \\
&\le1+\sum^n_{i=2}\frac{1}{i}\bigg[1+\frac{1}{1!}(\log i)^1+\cdots+\frac{1}{(k-1)!}
(\log i)^{k-1}\bigg] \\
&\le 1+\int^n_1\frac{1+\frac{1}{1!}(\log x)^1+\cdots+\frac{1}{(k-1)!}(\log x)^{k-1}}{x}dx.
\end{align*}
However
\[
\bigg[\frac{1}{x}\bigg(1+\frac{1}{1!}(\log x)^1+\cdots+\frac{1}{(k-1)!}(\log x)^{k-1}\bigg)\bigg]
=\bigg[\frac{1}{1!}(\log x)^1+\cdots+\frac{1}{k!}(\log x)^k\bigg]'.
\]
Thus, the inequality becomes:
\[
\le1+\int^n_1\bigg[\frac{1}{1!}(\log x)^1+\cdots+\frac{1}{k!}(\log x)^k\bigg]'dx=1+\frac{1}{1!}(\log n)^1+\cdots+\frac{1}{k!}(\log n)^k.
\]
Thus (\ref{sec2eq1}) becomes:
\[
T^{k+1}_{(n,m)}\le\frac{1}{n}\bigg[1+\frac{1}{1!}(\log n)^1+\cdots+\frac{1}{k!}(\log n)^k\bigg]
\]
which completes the induction argument.  \qb
\end{Proof}
\setcounter{cor}{3}
\begin{cor}\label{cor2.4}
For each $k,m\in\N^\ast$ $\dis\lim_{n\ra\infty}T^k_{(n,m)}=0$.
\end{cor}
\begin{Proof}
This is a direct consequence of Proposition \ref{prop2.3} since for each $k\in\N^\ast$
\[
\lim_{n\ra\infty}\frac{1}{n}\bigg[1+\frac{1}{1!}(\log n)^1+\cdots+\frac{1}{(k-1)!}(\log n)^{k-1}\bigg]=0.
\]
Indeed for each $\la\in\N^\ast$ we have
$\dis\lim_{n\ra\infty}\frac{(\log n)^\la}{n}=0$.  \qb
\end{Proof}
\setcounter{prop}{4}
\begin{prop}\label{prop2.5}
Let $k,n,\la\in\N^\ast$, then for each $i=1,2,\ld,\la$ it is true that:
\[
T^k_{(n+\la,n+i)}\ge\frac{(\la+1-i)^{k-1}}{(n+\la)^k}\cdot\frac{1}{(k-1)!}.
\]
\end{prop}
\begin{Proof}
For $k=1$, this means: $T_{(n+\la,n+i)}\ge\frac{1}{n+\la}$ and
equality holds for $i=1,\ldots,\la$.

The statement is true for $k=1$. We suppose it is true for
$k\in\N^\ast$ and we will show it is true for $k+1$:
\[
T^{k+1}_{(n+\la,n+i)}=\sum^\infty_{j=1}T_{(n+\la,j)}\cdot T^k_{(j,n+i)}
\]
and since both matrices are lower-triangular
\[
=\sum^{n+\la}_{j=n+i}T_{(n+\la,j)}\cdot T^k_{(j,n+i)}=\sum^{n+\la}_{j=n+i}
\frac{1}{n+\la}\cdot T^k_{(j,n+i)}.
\]
However, by the inductive argument, for $j=n+i,\ld,n+\la$
\[
T^k_{(j,n+i)}\ge\frac{(j-(n+i)+1)^{k-1}}{j^k}\cdot\frac{1}{(k-1)!}.
\]
Thus,
\begin{align*}
T^{k+1}_{(n+\la,n+i)}&=\frac{1}{n+\la}\cdot\sum^{n+\la}_{j=n+i}\frac{(j-(n+i)+1)^{k-1}}
{j^k}\cdot\frac{1}{(k-1)!} \\
&\ge\frac{1}{(n+\la)^{k+1}}\cdot\frac{1}{(k-1)!}\sum^{n+\la}_{j=n+i}(j-(n+i)+1)^{k-1} \ \ \text{since} \ \ j\le n+\la.
\end{align*}
Now, setting $\ti=j-(n+i)$ it becomes:
\[
=\frac{1}{(n+\la)^{k+1}}\cdot\frac{1}{(k-1)!}\cdot\sum^{\la-i+1}_{\ti=1}\ti^{k-1}.
\]
But
\[
\sum^{\la-i+1}_{\ti=1}\ti^{k-1}\ge\int^{\la-i+1}_0x^{k-1}dx=\frac{(\la-i+1)^k}{k}.
\]
Thus
\[
\ti^{k+1}_{(n+\la,n+i)}\ge\frac{1}{(n+\la)^{k+1}}\cdot\frac{1}{(k-1)!}\cdot\frac{(\la-i+1)^k}
{k}=\frac{(\la-i+1)^k}{(n+\la)^{k+1}}\cdot\frac{1}{k!}
\]
which completes the inductive argument.  \qb
\end{Proof}
\setcounter{cor}{5}
\begin{cor}\label{cor2.6}
Let $k,n,\la\in\N^\ast$. Then
\[
\sum^\la_{i=1}T^k_{(n+\la,n+i)}\ge\bigg(\frac{\la}{n+\la}\bigg)^k\cdot\frac{1}{k!}.
\]
\end{cor}
\begin{Proof}
By Proposition \ref{prop2.5}
\begin{align*}
\sum^\la_{i=1}T^k_{(n+\la,n+i)}&\ge\sum^\la_{i=1}\frac{(\la-i+1)^{k-1}}{(n+\la)^k}
\cdot\frac{1}{(k-1)!}\\
&=\frac{1}{(n+\la)^k}\cdot\frac{1}{(k-1)!}\cdot\sum^\la_{T=1}
T^{k-1}
\ge\frac{1}{(n+\la)^k}\cdot\frac{1}{(k-1)!}\cdot\int^\la_0x^{k-1}dx \\
&=\frac{1}{(n+\la)^k}\cdot\frac{1}{(k-1)!}\cdot\frac{\la^k}{k}=
\bigg(\frac{\la}{n+\la}\bigg)^k\cdot\frac{1}{k!}. \ \   \text{\qb}
\end{align*}
\end{Proof}
\setcounter{obs}{6}
\begin{obs}\label{obs2.7}
By Corollary \ref{cor2.2} if $k,n,\la\in\N^\ast$ and $\la\le n$,
then for each $i=1,\ld,\la$ it is true that
$T^k_{(n+\la,n+i)}\le\frac{2}{n+\la}$.

Thus, $\ssum^\la_{i=1}T^k_{(n+\la,n+i)}\le2\cdot\frac{\la}{n+\la}$.
\end{obs}
\setcounter{prop}{7}
\begin{prop}\label{prop2.8}
Let $k,n,\la,\La\in\N^\ast$, $\la<\La$ and $\La\le n$. Then for $\ti=1,2,\ld,\la$ $a_\ti=\frac{T^k_{(n+\la,n+\ti)}}{T^k_{(n+\La,n+\ti)}}$ is decreasing with respect to $\ti$. Moreover it is true that $2\ge a_1\ge\cdots\ge a_\la>0$.
\end{prop}
\begin{Proof}
In order to prove that $a_\ti$ is decreasing, the case $\La=\la+1$ (for arbitrary $\la$ such that $\la+1\le n)$ is enough.

Indeed
\[
a_\ti=\frac{T^k_{(n+\la,n+\ti)}}{T^k_{(n+\la+1,n+\ti)}}\cdot\frac{T^k_{(n+\la+1,
n+\ti)}}{T^k_{(n+\la+2,n+\ti)}}\cdot\,\cdots\,\cdot\frac{T^k_{(n+\La-1,n+\ti)}}
{T^k_{(n+\La,n+\ti)}}
\]
would be a product of decreasing sequences (that have finitely many
terms) of positive real numbers, thus $a_\ti$ would be decreasing.

To prove the case $\La=\la+1$, we proceed by induction on $k$. For $k=1$:
\[
a_\ti=\frac{T_{(n+\la,n+\ti)}}{T_{(n+\la+1,n+\ti)}}=\frac{\frac{1}{n+\la}}
{\frac{1}{n+\la+1}}=\frac{n+\la+1}{n+\la}
\]
which is constant with respect to $\ti$.

Thus, it is true that $a_1\ge a_2\ge\cdots\ge a_\la$ for $k=1$. (In
fact, equality holds for $k=1$).

We assume that $a_\ti$ is decreasing for $k\in\N^\ast$ and we will
show that it is decreasing for $k+1$.
$T^{k+1}_{(n+\la,n+\ti)}=\ssum^\infty_{i=1}T_{(n+\la,i)}\cdot
T^k_{(i,n+\ti)}$ but $T$, $T^k$ are lower-triangular, thus:
\[
=\sum^{n+\la}_{i=n+\ti}T_{(n+\la,i)}\cdot T^k_{(i,n+\ti)}=\frac{1}{n+\la}\cdot\Big[T^k_{(n+\ti,n+\ti)}+\cdots+T^k_{(n+\la,n+\ti)}\Big].
\]
Thus $a_\ti$, in the case $(k+1)$ can be written as:
\begin{eqnarray}
a_\ti=\frac{T^{k+1}_{(n+\la,n+\ti)}}{T^{k+1}_{(n+\la+1,n+\ti)}}=
\frac{\frac{1}{n+\la}\Big[T^k_{(n+\ti,n+\ti)}+\cdots+T^k_{(n+\la,n+\ti)}\Big]}
{\frac{1}{n+\la+1}\Big[T^k_{(n+\ti,n+\ti)}+\cdots+T^k_{(n+\la+1,n+\ti)}\Big]} \label{sec2eq2}
\end{eqnarray}
and we need to show that $a_\ti\ge a_{\ti+1}$ $(t=1,\ld,\la-1)$. That means:
\[
\frac{\frac{1}{n+\la}\Big[T^k_{(n+\ti,n+\ti)}\!+\!\cdots\!+\!T^k_{n+\la,n+\ti)}\Big]}
{\frac{1}{n+\la+1}\Big[T^k_{(n+\ti,n+\ti)}\!+\!\cdots\!+\!T^k_{(n+\la+1,n+\ti)}\Big]}\!\ge\!
\frac{\frac{1}{n+\la}}{\frac{1}{n+\la+1}}
\frac{\Big[T^k_{(n+\ti+1,n+\ti+1)}\!+\!
\cdots\!+\!T^k_{(n+\la,n+\ti+1)}\Big]}{\Big[T^k_{(n+\ti+1,n+\ti+1)}\!+\!\cdots\!+\!T^k_{(n+\la+1,n+\ti+1)}\Big]}.
\]
Now if we set $A=T^k_{(n+\ti,n+\ti)}+\cdots+T^k_{(n+\la,n+\ti)}$ and $B=T^k_{(n+\ti+1,n+\ti+1)}+\cdots+T^k_{(n+\la,n+\ti+1)}$ the previous expression is written as:
\begin{align*}
&A(B+T^k_{(n+\la+1,n+\ti+1)})\ge B(A+T^k_{(n+\la+1,n+\ti)})\,\Leftrightarrow \\
&A\cdot T^k_{(n+\la+1,n+\ti+1)}\ge B\cdot T^k_{(n+\la+1,n+\ti)}
\end{align*}
which is true by the induction argument since for $\la\ge\rho\ge\ti+1$
\begin{align*}
&T^k_{(n+\rho,n+\ti)}\cdot T^k_{(n+\la+1,n+\ti+1)}\ge T^k_{(n+\rho,n+\ti+1)}\cdot T^k_{(n+\la+1,n+\ti)}\, \Leftrightarrow \\
&\frac{T^k_{(n+\rho,n+\ti)}}{T^k_{(n+\la+1,n+\ti)}}\ge\frac{T^k_{(n+\rho,n+\ti+1)}}
{T^k_{(n+\la+1,n+\ti+1)}}
\end{align*}
is exactly the inductive argument for $k$.

Moreover, we needn't take into account the part
$T^k_{(n+\ti,n+\ti)}$ of $A$.

Finally, the expression (\ref{sec2eq2}) for $a_\ti$ implies that $a_\ti<\frac{n+\la+1}{n+\la}$ since the sum in the denominator is bigger than the sum in the nominator.

Thus, in the general case for $\la,\La$ it will be true that
\[
a_\ti<\frac{n+\la+1}{n+\la}\cdot\frac{n+\la+2}{n+\la+1}\cdots\frac{n+\La}{n+\La-1}=
\frac{n+\La}{n+\la}<2
\]
since $\La\le n$ and $\La,\la>0$.

It is also obvious, by definition, that $a_\ti>0$ $\ti=1,2,\ld,\la$. Thus, indeed, $2\ge a_1\ge a_2\ge\cdots\ge a_\la>0$.  \qb
\end{Proof}
\begin{prop}\label{prop2.9}
Let $V$ be a vector space over $\R$ and $\|\,\|$ a semi-norm on $V$.
Let $\la\in\N$ and $\bi_1,\ld,\bi_\la\in V$ such that there exist
$M\in\R_+$: such that for all $j=1,2,\ld,\la$ we have
$\Big\|\ssum^j_{\ti=1}\bi_\ti\Big\|\le M$.
\end{prop}

If $a_1,\ld,a_\la\in\R_+$, $a_1\ge a_2\ge\cdots\ge a_\la>0$ then for
every $j=1,2,\ld,\la$ we have
$\Big\|\ssum^j_{\ti=1}a_\ti\bi_\ti\Big\|\le a_1\cdot M$.
\begin{Proof}
By Abel's transformation
\[
\sum^j_{\ti=1}a_\ti\bi_\ti=a_j(\bi_1+\cdots+\bi_j)-\sum^{j-1}_{\ti=1}(\bi_1+\cdots+
\bi_\ti)(a_{\ti+1}-a_\ti).
\]
Considering the seminorm we get:
\begin{align*}
\bigg\|\sum^j_{\ti=1}a_\ti\bi_\ti\bigg\|&\le a_j\|\bi_i+\cdots+\bi_j\|+\sum^{j-1}_{\ti=1}\|\bi_1+\cdots+\bi_\ti\|(a_\ti-a_{\ti+1}) \\
&\le a_j\cdot M+\sum^{j-1}_{\ti=1}M\cdot(a_\ti-a_{\ti+1})=a_1M. \ \ \text{\qb}
\end{align*}
\end{Proof}

Now, let $\vthi\in V^\N$ be a sequence $(\vthi=(\vthi_n))$. We prove
the following:
\begin{prop}\label{prop2.10}
Let $n,k,\la\in\N^\ast$. Then
\[
[T^k(\vthi)]_{n+a}=\frac{n}{n+a}\cdot[T^k(\vthi)]_n+
\frac{[T^{k-1}(\vthi)]_{n+1}+\cdots+[T^{k-1}(\vthi)]_{n+a}}{n+a}.
\]
\end{prop}
\begin{Proof}
Indeed,
\begin{align*}
[T^k(\vthi)]_{n+a}&=\frac{[T^{k-1}(\vthi)]_1+\cdots+[T^{k-1}(\vthi)]_{n+a}}{n+a} \\
&=\frac{\frac{n}{n}\Big[[T^{k-1}(\vthi)]_1+\cdots+[T^{k-1}(\vthi)]_n\Big]+
[T^{k-1}(\vthi)]_{n+1}+\cdots+[T^{k-1}(\vthi)]_{n+a}}{n+a} \\
&=\frac{n\cdot[T^k(\thi)]_n+[T^{k-1}(\vthi)]_{n+1}+\cdots+[T^{k-1}(\vthi)]_{n+a}}
{n+a}
\end{align*}
which concludes the proof. \qb
\end{Proof}
\setcounter{cor}{10}
\begin{cor}\label{cor2.11}
Let $n,k,\la\in\N^\ast$.\ Then $[T^k(\vthi)]_{n+a}$ is a convex combination of $[T^k(\vthi)]_n,$ $[T^{k-1}(\vthi)]_n,\ld,[T(\vthi)]_n,\vthi_{n+1},\ld,\vthi_{n+a}$.
\end{cor}
\begin{Proof}
Indeed,by Proposition \ref{prop2.10} $[T^k(\vthi)]_{n+a}$ is a
convex combination of $[T^k(\vthi)]_n$,
$[T^{k-1}(\vthi)]_{n+1},\ld,[T^{k-1}(\vthi)]_{n+a}$.

However, each of the $[T^{k-1}(\vthi)]_{n+i}$ $(i=1,\ld,a)$ is a convex combination of\\ $[T^{k-1}(\vthi)]_n$, $[T^{k-2}(\vthi)]_{n+1},\ld,[T^{k-2}(\vthi)]_{n+i}$.

Thus, repeating this step $k$ times we see that $[T^k(\vthi)]_{n+a}$
is a convex combination of $[T^k(\vthi)]_n$,
$[T^{k-1}(\vthi)]_n,\ld,[T(\vthi)]_n$ and
$\vthi_{n+1},\ld,\vthi_{n+a}$. \qb
\end{Proof}

Finally, we mention the following
\setcounter{rem}{11}
\begin{rem}\label{rem2.12}
If $V$ is a locally convex vector space and $a$ is a convergent
sequence in $V$ towards a limit $\el\in V$, then $[T(a)]_n\ra\el$,
as $n\ra+\infty$, as well.
\end{rem}
\section{A generic result on the iterates of the Cesaro operator} \label{sec3}
\setcounter{Def}{-1}
\begin{Def}\label{Def3.0}
Let $V$ be a metrizable, separable, locally convex vector space over
$\R$. The topology on $V$ is given by a countable family of
separated semi-norms. It can be generated by the metric
$d(x,y)=\ssum^\infty_{i=1}\frac{1}{2^i}\frac{\|x-y\|_i}{1+\|x-y\|_i}$
$(x,y\in V$, $\|\,\|_i: i=1,2,\ld$ are the semi-norms).

For each $\e>0$ we define $N_\e$ to be the smallest natural number
such that $\ssum^\infty_{i=N_\e}\frac{1}{2^i}<\e$. The number $N_\e$
is in this way the number of $\e$-important seminorms.

For instance, if $x,y\in V$ and $\|x-y\|_i\le\e$ $i=1,2,\ld,N_\e$
then $d(x,y)<\e+\e=2\e$.

Moreover, let $T:V^\N\ra V^\N$ be the Cesaro operator: if
$a=(a_n)\in V^\N$ then $[T(a)]_n=\frac{a_1+\cdots+a_n}{n}$. Also
$T^k=T\circ T^{k-1}$ $k=2,3,\ld\;.$
\end{Def}

I prove the following:
\setcounter{thm}{0}
\begin{thm}\label{thm3.1}
If $V$ and $T$ are defined as in Definition \ref{Def3.0}. and if
$A\subseteq V$ is an arbitrary subset, then the set:
\[
\Thi_{A,V}=\left\{\begin{array}{l}
                    \vthi\in A^\N:\,[T^k(\vthi)]_n \ \ \text{is dense in} \\
                    \text{conv}(A) \ \ \fa\;k=1,2,\ld
                  \end{array}\right\}
\]
is dense and $G_\de$ in $A^\N$ under the product topology.
\end{thm}

We will first prove a simpler case:
\begin{thm}\label{thm3.2}
Let $W$ be a separable Frechet space (completely metrizable, locally
convex) and $B\subset W$ a closed subset of $W$. Then the set:
\[
\Thi_{B,W}=\left\{\begin{array}{l}
                    \vthi\in B^\N:\,[T^k(\vthi)]_n \ \ \text{is dense in} \\
                    \text{conv}(B) \ \ \fa\;k=1,2,\ld
                  \end{array}\right\}
\]
is dense and $G_\de$ in $B^\N$ under the product topology.
\end{thm}
\noindent
{\bf Proof of Theorem \ref{thm3.2}.} Since $W$ is separable and
metrizable and conv$B\subseteq W$, then conv$B$ will also be
separable. Let $\big\{x_\la\big\}^\infty_{\la=1}$ be a sequence that
is dense in conv$(B)$ and $x_\la\in$conv$(B)$ for all
$\la=1,2,\ld\;.$

For each $k,\la\in\N$ I define the set $\Thi_{x_\la,k}$ as follows:
\[
\Thi_{x_\la,k}=\left\{\begin{array}{l}
                        \vthi\in B^\N:\,\exists\;n_0\in\N \ \ \text{such} \\
                        \text{that} \ \ d([T^k(\vthi)]_{n_0},x_\la)<\frac{1}{\la}
                      \end{array}\right\}.
\]
It is true that $\Thi_{B,W}$ can be expressed as:
\[
\Thi_{B,W}=\bigcap^\infty_{k=1}\bigcap^\infty_{\la=1}\Thi_{x_\la,k}.
\]
This is an obvious consequence of the definition of $\Thi_{B,W}$.

Moreover,
\[
\Thi_{x_\la,k}=\bigcup^\infty_{m=1}(T^k)^{-1}\bigg(W\times W\times\cdots\times W\times \underbrace{B\bigg(x_\la,\frac{1}{\la}\bigg)}_{m-\text{position}}\times W\times\cdots\bigg)
\]
where $T^k:B^\N\ra W^\N$ is continuous and
$B\Big(x_\la,\frac{1}{\la}\Big)$ is the open ball of center $x_\la$
and radius $\frac{1}{\la}$.

The set $W\times W\times\cdots\times W\times
B\Big(x_\la,\frac{1}{\la}\Big)\times W\times\cdots$ is open in
$W^\N$ under the product topology, thus its inverse image under the
continuous map $T^k$ will also be open.

Thus $\Thi_{x_\la,k}$ is open as a union of open
sets.\vspace*{0.2cm}

We will now show that $\Thi_{x_\la,k}$ is also dense in $B^\N$. This
is a direct consequence of the following lemma.
\setcounter{lem}{2}
\begin{lem}\label{lem3.3}
Let $W$ as in Theorem \ref{thm3.2} and $A\subseteq W$ an arbitrary
subset, $\vthi_1,\ld,\vthi_\rho\in A$ $\rho$-many given terms
(without necessarily having $\rho\neq0$), $x\in$conv$A$ ``target'',
$\e>0$ and $k\in\N$ be given.

Then, there exist $n>\rho$ $(n\in\N)$ and
$\vthi_{\rho+1},\ld,\vthi_n\in A$ terms, such that if $\vthi'$ is a
sequence, $\vthi'\in A^\N$ that ``starts with''
$\vthi_1,\ld,\vthi_n$ (that means $\vthi'_i=\vthi_i$, $i=1,\ld,n$)
it is true that $d([T^k(\vthi)]_n,x)<\e$.

We mention that $W$ need not be completely metrizable for the lemma to hold.
\end{lem}
\noindent
{\bf Proof of Lemma \ref{lem3.3}.} I will first prove that there is
a sequence $\vthi'\in A^\N$ that ``starts with''
$\vthi_1,\ld,\vthi_\rho$ and has a limit in conv$A$ as $n\ra\infty$.
(That is: $\dis\lim_{n\ra\infty}[T(\vthi')]_n$ also exists and
belongs to conv$(A)$) and moreover, that
$d\Big(\dis\lim_{n\ra\infty}[T(\vthi')]_n,x\Big)<\frac{2\cdot\e}{3}$.

$x\in$conv$(A)$ thus $\exists\;v\in\N$ and $\la_1,\ld,\la_v\in\R$
and $A_1,\ld,A_v\in A$ such that $\la_i>0$, $\ssum^v_{i=1}\la_i=1$
and $x=\ssum^v_{i=1}\la_iA_i$.

For each $m\in\N$, there exist $m_1,\ld,m_v\in\N$ such that
$\ssum^v_{i=1}m_i=m$ and
$\Big|\la_i-\frac{m_i}{m}\Big|\le\frac{1}{m}$.

Indeed, to show this, for each $m$, I consider maximal natural
numbers $m'_i$ such that $\frac{m'_i}{m_i}<\la_i$. Furthermore, I
define $m_i$ to be either $m'_i$ or $m'_i+1$ in such a way that
$\ssum^v_{i=1}m_i=m$.

Since this is true for each $m\in\N$, I can choose $m$ to be such
that $\frac{\|A_i\|_j}{m}<\frac{\e}{3v}$, $i=1,\ld,v$,
$j=1,2,\ld,N_{\e/3}$. Then
\[
\bigg\|x-\sum^v_{i=1}\frac{m_i}{m}\cdot A_i\bigg\|_j\!=\!\bigg\|\sum^v_{i=1}A_i\bigg(\la_i-\frac{m_i}{m}\bigg)\bigg\||_j \le\frac{1}{m}\cdot\sum^v_{i=1}\|a_i\|_j<v\cdot\frac{\e}{3v}=\frac{\e}{3}, \ \ j=1,2,\ld,N_{\e/3}.
\]
I define $x'=\ssum^v_{i=1}\frac{m_i}{m}A_i$. Now $\|x-x'\|_j<\e/3$
$j=1,\ld,N_{\e/3}$ thus $d(x,x')<\frac{2\cdot\e}{3}$.

I consider the sequence $\vthi'\in A^\N$ defined as follows: $\vthi'_i=\vthi_i$ $i=1,2,\ld,\rho$ and from index $(\rho+1)$ onwards I repeat the $m$-tuple of elements of $A$
\[
\underbrace{\overbrace{A_1A_1\cdots A_1}^{m_1\text{-many}}\overbrace{A_2\cdots A_2}^{m_2\text{-many}}\cdots A_{v-1}\overbrace{A_v\cdots A_v}^{m_v\text{-many}}}_{m\text{-many}}\qquad (m_1+\cdots+m_v=m).
\]

I will show that $\dis\lim_{n\ra\infty}[T(\vthi')]_n=x'$.

Let $u<m$ be a natural number and $q$ an arbitrary natural number. I
examine the term of index $(\rho+qm+u)$ of $T(\vthi')$:
\[
[T(\vthi')]_{\rho+qm+u}=\frac{\vthi_1+\cdots+\vthi_\rho}{\rho+qm+u}+
\frac{q\Big(\ssum^v_{i=1}m_iA_i\Big)}{\rho+qm+u}+\frac{\overbrace{A_1+\cdots A_\si}^{u\text{-many}}}{\rho+qm+u}
\]
where $\si$ is an index $\si\le v$. Moreover
$\ssum^v_{i=1}m_iA_i=m\cdot x'$. Thus for every $j\in\N$ it is true
that:
\[
\|[T(\vthi')]_{\rho+qm+u}-x'\|_j\le\bigg\|\frac{\vthi_1+\cdots+\vthi_\rho}{\rho+qm+u}\bigg\|_j+
\bigg\|\frac{x'\cdot(\rho+u)}{\rho+qm+u}\bigg\|_j+\bigg\|\frac{A_1+\cdots+A_\si}{\rho+qm+u}\bigg\|_j.
\]
For each $j\in\N$, the real numbers
$\|A_1\|_j,\|A_1+A_1\|_j,\ld,\|\underbrace{A_1+\cdots+A_v}_{m\text{-many}}\|_j$
are finitely many ($m$-many) positive real numbers. Thus, one of
them is maximal. The previous inequality is also true for each
$u=1,2,\ld,m$ and each $q\in\N$. Moreover, the right part approaches
0 when $q\ra\infty$, independent of the value of $u$. Finally since
there are finitely many options for $u$, we deduce that
\[
\lim_{n\ra\infty}\|[T(\vthi')]_n-x'\|_j=0 \qquad j=1,2,\ld\;.
\]
Since this is true for each $j\in\N$, the following also holds:
$\dis\lim_{n\ra\infty}d([T(\vthi')]_n,x')=0$.

The Cesaro operator does not affect the limit of a convergent
sequence, thus for each $i=1,2,\ld$ it will be true that
$\dis\lim_{n\ra\infty}d([T^i(\vthi')]_n,x')=0$. Thus, it will be
true for $i=k$. Thus, there exists a natural number $n_0$ such that
$d([T^k(\vthi')]_{n_0},x')<\frac{\e}{3}$. Now, we also showed that
$d(x,x')<\frac{2\e}{3}$.

Thus, $d([T^k(\vthi')]_{n_0},x)<\e$.

The number $n_0$ and the terms $\vthi'_{\rho+1},\ld,\vthi'_{n_0}$ satisfy the argument and this concludes the proof of Lemma \ref{lem3.3}.  \qb\vspace*{0.2cm}

I continue now, with the Proof of Theorem \ref{thm3.2}.\vspace*{0.2cm}

Using Lemma \ref{lem3.3} it is easy now, to deduce that
$\Thi_{x_\la,k}$ is dense in $B^\N$:

Let $\bi\in B^\N$ be an arbitrary sequence and let $\rho\in\N$. I
consider Lemma \ref{lem3.3}, setting $\vthi_i=\bi_i$,
$i=1,\ld,\rho$, $x=x_\la$, $\e=\frac{1}{\la}$ and $k$ to be the one
that defines $\Thi_{x_\la,k}$. Now, I consider an arbitrary sequence
$\vthi'\in B^\N$ that ``starts with'' the terms
$\vthi_1,\ld,\vthi_n$ that are provided by Lemma \ref{lem3.3}.

Then $\vthi'\in\Thi_{x_\la,k}$ (Indeed, the index $n_0$ that is
demanded by the definition of $\Thi_{x_\la,k}$ is the number $n$
provided by Lemma \ref{lem3.3}).

Also by definition $\vthi'_i=\bi_i$ $i=1,\ld,\rho$ and since such a
sequence $\vthi'$ can be constructed for an arbitrary $\bi\in B^\N$
and $\fa\;\rho\in\N$, $\Thi_{x_\la,k}$ is dense in $B^\N$.

I have shown that
$\Thi_{B,W}=\bigcap\limits^\infty_{k=1}\bigcap\limits^\infty_{\la=1}\Thi_{x_\la,k}$
and proved that $\Thi_{x_\la,k}$ is open and dense in $B^\N$ under
the product topology. Since $B$ is closed in the complete space $W$
it will also be complete under the metric $d$. Thus, $B^\N$ is
completely metrizable, that means it is a Baire space.

Thus from Baire's Category Theorem I deduce that $\Thi_{B,W}$ is dense and $G_\de$ in $B^\N$. This completes the proof of Theorem \ref{thm3.2}.  \qb \vspace*{0.2cm}

I proceed now with the proof of Theorem \ref{thm3.1}. \vspace*{0.2cm} \\
\noindent
{\bf Proof of Theorem \ref{thm3.1}.} I consider a completion $W$ of $V$ and the closure $\bar{A}^W$ of $A$ in $W$. Then because of Theorem \ref{thm3.2} the set
\[
\Thi_{\bar{A}^W,W}=\left\{\begin{array}{l}
                            \vthi\in(\bar{A}^W)^\N:\,[T^k(\vthi)]_n \ \ \text{is dense} \\
                            \text{conv}(\bar{A}^W) \ \ \fa\;k=1,2,\ld
                          \end{array}\right\}
\]
is dense and $G_\de$ in $(\bar{A}^W)^\N$.

It is true that
\setcounter{equation}{0}
\begin{eqnarray}
\Thi_{\bar{A}^W,W}\cap A^\N=\Thi_{A,V}.  \label{sec3eq1}
\end{eqnarray}
Indeed, the direction $\subseteq$ is obvious since, if $\vthi\in A^\N$, then $T^k(\vthi)\in(\text{conv}(A))^\N$ $\fa\;k=1,2,\ld$ and if $[T^k(\vthi)]_n$ is dense in conv$(\bar{A}^W)$ it will also be dense in conv$(A)$ (since conv$(A)\subseteq\,$conv$(\bar{A}^W)$).

For the other direction $\supseteq$ we have:

If $\vthi\in\Thi_{A,V}$ then $\vthi\in A^\N$  by definition. Now,
$[T^k(\vthi)]_n$ is dense in conv$(A)$ for every $k=1,2,\ld$ and
conv$(A)$ is dense in conv$(\bar{A}^W)$ since $A$ is dense in
$\bar{A}^W$, thus $[T^k(\vthi)]_n$ is dense in conv$(\bar{A}^W)$,
thus $\vthi\in\Thi_{\bar{A}^W,W}$.

(\ref{sec3eq1}) implies that $\Thi_{A,V}$ is $G_\de$ in $A^\N$, because $\Thi_{\bar{A}^W,W}$ is $G_\de$ in $(\bar{A}^W)^\N$.

We will now show that $\Thi_{A,V}$ is dense in $A^\N$. It suffices
to show that $\Thi_{\bar{A}^W,W}\cap A^\N$ is dense in
$\Thi_{\bar{A}^W,W}$ (which is dense in $(\bar{A}^W)^\N\supseteq
A^\N)$.

Indeed let $\vthi'\in\Thi_{\bar{A}^W,W}$.

Then, there exist infinitely many sequences $\vthi^i=(\vthi^i)_n$
and $\si^i=(\si^i)_n$ $i=1,2,\ld$ such that
$(\vthi')_n=(\vthi^i)_n+(\si^i)_n$ and
$\dis\lim_{n\ra\infty}(\si^i)_n=0$ and
$\dis\lim_{i\ra\infty}(\si^i)_v=0$, $\fa\;v=1,2,\ld\;.$

Since $[T^k(\vthi')]_n=[T^k(\vthi^i)]_n+[T^k(\si^i)]_n$ and since
$[T^k(\si^i)]_n\ra0$. (The Cesaro operator does not affect the limit
of convergent sequences), then $[T^k(\vthi^i)]_n$ will also be dense
in conv$(\bar{A}^W)$. That is $\vthi^i\in\Thi_{\bar{A}^W,W}$,
$\fa\;i=1,2,\;d\;.$

Since $\dis\lim_{i\ra\infty}(\si^i)_v=0$ $\fa\;v=1,2,\ld$ and since
$(\vthi')_n=(\vthi^i)_n+(\si^i)_n$, it is true that
$\dis\lim_{i\ra\infty}\vthi^i=\vthi'$ where the limit is taken in
$(\bar{A}^W)^\N$ under the product topology. But $\vthi^i\in A^\N$
$\fa\;i=1,2,\ld\;.$ That means that $\vthi'$ belongs to the closure
of $A^\N$ in $(\bar{A}^W)^\N$. Since this is true for an arbitrary
$\vthi'\in\Thi_{\bar{A}^W,W}$ I deduce that $\Thi_{\bar{A}^W,W}\cap
A^\N=\Thi_{A,V}$ is dense in $\Thi_{\bar{A}^W,W}$ thus, also dense
in $A^\N$.

This completes the proof of Theorem \ref{thm3.1}.  \qb
\setcounter{rem}{3}
\begin{rem}\label{rem3.4}
The fact  that $\Thi_{B,W}$ is $G_\de$-dense is important because it
lies in $B^\N$ which is a complete space. In a complete space, a
$G_\de$-dense set is quite ``big''.

However,  although the general theorem gives us more information,
the extra knowledge is not significant since in a not complete
$(A^\N)$ space, the $G_\de$-dense sets need not be big. Thus the
general theorem could be seen as an existence theorem, i.e.
$\Thi_{A,V}\neq\emptyset$.
\end{rem}
\section{Main Result}\label{sec4}
\noindent

Let $V$, $N_\e$, $T$ be as in Definition \ref{Def3.0}.
\setcounter{Def}{-1}
\begin{Def}\label{Def4.0}
I consider the function
$\de:\Big(0,\frac{1}{2}\Big)\ra\Big(0,\frac{1}{2}\Big)$ such that if
$x\in B(0,\de(\e))$ (the open ball with center 0 and radius
$\de(\e)$ in $V$) then $\|x\|_\rho<\e$, $\rho=1,2,\ld,N_\e$. Such a
function clearly exists and one suitable choice would be
$\de(\e)=\frac{\e}{2^{N_\e+1}}$.
\end{Def}

I will prove the following:
\setcounter{thm}{0}
\begin{thm}\label{thm4.1}
Let $(V,d)$ be a separable, locally convex metric vector space, $\mathfrak{M}\subseteq\N$ an infinite set of natural numbers, and $A\subset V$ a subset such that $\overline{\text{conv}A}=V$.

Then the set
\[
\Thi(A,\mf)=\left\{\begin{array}{l}
                        (\vthi_n)\in A^\N:\;\text{for each}\; (x_k)\in V^\N \;\text{sequence of ``targets'', there exists} \\
                        \text{a sequence of indices}\;\la_n, n=1,2,\ld,\;\text{such that}\; \la_n\in\mf\;\text{and}\\
                        \lim[T^k(\vthi)]_{\la_n}=x_k \; \fa\;k=1,2,\ld
                      \end{array}\right.
\]
is dense and $G_\de$ in $A^\N$. Moreover, if $A$ is a vector space,
then $\Thi(A,\mf)$ contains a linear subspace of $V^\N$ with the
exception of 0.
\end{thm}

This is a consequence of the following finite version.
\begin{thm}\label{thm4.2}
Let $V,A,\mf$ be as in Theorem \ref{thm4.1}. Let $k\in\N$ and
$\vthi_1,\ld,\vthi_\rho\in A$ $\rho$-many terms (without necessarily
having $\rho\neq0$), a $k$-tuple of ``targets'', $x_1,\ld,x_k\in V$
and $\e>0$ be given. Then there exists $n\in\mf$, $n>\rho$ and
$\vthi_{\rho+1},\ld,\vthi_n\in A$ terms such that if a sequence
$\vthi\in V^\N$ ``starts with'' the terms $\vthi_1,\ld,\vthi_n$ then
it is true that $d([T^i(\vthi)]_n,x_i)<\e$ $i=1,\ld,k$.
\end{thm}

Supposing Theorem \ref{thm4.2} is true I will proceed with the proof of Theorem \ref{thm4.1}. \vspace*{0.2cm}\\
\noindent
{\bf Proof of Theorem \ref{thm4.1}.} I will first prove that $\Thi(A,\mf)$ is nonempty.

Since $V$ is separable, $V^\N$ is also separable. Indeed if $\xi$ is
a dense countable subset of $V$, then the set of sequences of
elements of $\xi$ with finite support is dense in $V^\N$ under the
product topology.

Let $\big\{\X_\la\big\}^\infty_{\la=1}$ be a sequence of elements of
$V^\N$ that is dense in $V^\N$. This means that
$\X_\la\equiv(\X_\la)_n\in V^\N$ and the previous comment implies we
can further demand that $\fa\;\la=1,2,\ld$ it is true that
$\sharp\,\{n\in\N:(\X_\la)_n\neq0\}<\infty$. We can also demand that
$\X_\la$ is of the form $(x_1,\ld,x_k,0,0,\ld)$ for some $k\in\N$
and $x_1,\ld,x_k\in V^\ast$.

I will construct an element of $\Thi(A,\mf)$.\vspace*{0.2cm}

I consider an empty set of  initial terms $(\rho=0)$,
$\e=\frac{1}{1}$, ``targets'' $x_1,\ld,x_{k_1}$ to be the nonzero
terms of $\X_1$ and $k=k_1$. According to Theorem \ref{thm4.2} there
exists $n_1\in\mf$ and $\vthi_1,\ld,\vthi_n\in A$ such that if
$\vthi'\in V^\N$ ``starts with'' $\vthi_1,\ld,\vthi_{n_1}$ then it
is true that $d([T^i(\vthi')]_{n_1},x_i)<\frac{1}{1}$,
$i=1,2,\ld,k_1$.

For each  $\la=1,2,3,\ld$ I inductively define new terms as follows:
I consider $\vthi_1,\ld,\vthi_{n_\la}$ as the initial terms,
$x'_1,\ld,x'_{k_{\la+1}}$ the nonzero terms of $\X_{\la+1}$ as the
``targets'', $k=k_{\la+1}$, $\e=\frac{1}{\la+1}$. Now according to
Theorem \ref{thm4.2} there exist a natural number $n_{\la+1}\in\mf$,
$n_{\la+1}>n_\la$ and new terms
$\vthi_{n_\la+1},\ld,\vthi_{n_{\la+1}}\in A$ such that if $\vthi'\in
V^\N$ ``starts with'' $\vthi_1,\ld,\vthi_{n_{\la+1}}$, then it is
true that
\[
d([T^i(\vthi')]_{n_{\la+1}},x'_i)<\frac{1}{\la+1} \qquad i=1,\ld,k_{\la+1}.
\]

I consider the sequence $\vthi$ that is formed by all these terms $\vthi_i$ that I inductively defined. I will show that $\vthi\in\Thi(A,\mf)$.

It is clear that $\vthi_n\in A$ $\fa\;n=1,2,\ld\;.$ Now, let
$(x_k)\in V^\N$. Since $\X_n$ is dense in $V^\N$, there exist a
strictly increasing sequence of indices $\rho_v$ such that
$\dis\lim_{v\ra\infty}\X_{\rho_v}=(x_k)$. Because of the way $\vthi$
is constructed it is true that $d([T^i(\vthi)]_{n_{\rho_v}}$,
$(\X_{\rho_v})_i)<\frac{1}{\rho_v}$ where $i$ runs from 1 to the
last nonzero term of $\X_{\rho_v}$.

Let $k\in\N$
$\dis\lim_{v\ra\infty}d([T^k(\vthi)]_{n_{\rho_v}},x_k)=\dis\lim_{v\ra\infty}
d([T^k(\vthi)]_{n_{\rho_v}}$,
$(\X_{\rho_v})_k)\le\dis\lim_{v\ra\infty}\frac{1}{\rho_v}=0$, since
$x_k=\dis\lim_{v\ra\infty}(\X_{\rho_v})_k$ $\fa\;k=1,2,\ld\;.$

Since $k$ was arbitrary it will be true that
$\dis\lim_{v\ra\infty}d([T^k(\vthi)]_{n_{\rho_v}},x_k)=0$.
$\fa\;k=1,2,\ld\;.$

Since $(x_k)$ was also an arbitrary element of $V^\N$, it will be
true for all $(x_k)\in V^\N$  that there exists a strictly
increasing sequence of indices $\rho_v$ such that
$\dis\lim_{\rho_v\ra\infty} d([T^k(\vthi)]_{n_{\rho_v}},x_k)=0$.
Furthermore, $n_{\rho_v}\in\mf$ by the construction of $\vthi$ and
hence the set $\Thi(A,\mf)$ is nonempty indeed.

Since  $\Thi(A,\mf)\neq\emptyset$, we may consider an element
$\vthi$ of it. I will show that $\Thi(A,\mf)$ is dense in $A^\N$,
under the product topology. It is enough to show that if
$a\equiv(a_n)\in A^\N$ is an arbitrary sequence, then for each
$v\in\N$, there exists a sequence $\vthi'\in\Thi(A,\mf)$ such that
$\vthi'_i=a_i$ $i=1,\ld,v$.

I consider the given sequence $(a)$ and an arbitrary number $v\in\N$ and I define $\vthi'$ to be: $\vthi'_i=a_i$ $i=1,\ld,v$, $\vthi'_i=\vthi_i$ $i\ge v+1$.

I will show that $\vthi'\in\Thi(A,\mf)$.

Clearly $(\vthi')_n\in A$ $\fa\;n$.

Let $(x_k)\in V^\N$ be a sequence of targets and $\la_n$ the
sequence of indices for which
$\dis\lim_{n\ra\infty}d([T^k(\vthi)]_{\la_n},\;x_k)=0$.

Let $k,j\in\N^\ast$ be arbitrary natural numbers. Then
\[
\|[T^k(\vthi)]_{\la_n}-[T^k(\vthi')]_{\la_n}\|_j\le\sum^v_{i=1}T^k_{(\la_n,i)}\cdot
\|\vthi_i\|_j+\sum^v_{i=1}T^k_{(\la_n,i)}\|a_i\|_j
\]
where $T^k_{(n,i)}$ is the $(n,i)$ element of the matrix of $T^k$.

The second part of the inequality tends to 0, as $n$ approaches infinity since $\dis\lim_{n\ra\infty}T^k_{(n,i)}=0$ $\fa\;i=1,2,\ld$ by Corollary \ref{cor2.4} and since both sums are finite. Thus,
\[
\lim_{n\ra\infty}\|[T^k(\vthi)]_{\la_n}-[T^k(\vthi')]_{\la_n}\|_j=0 \ \ \fa\;k,j\in\N.
\]
Thus
\[
\lim_{n\ra\infty}d([T^k(\vthi)]_{\la_n},\;[T^k(\vthi')]_{\la_n})=0 \ \ \fa\;k\in\N.
\]
Thus
\[
\lim_{n\ra\infty}d([T^k(\vthi')]_{\la_n},\;x_n)=0 \ \ \fa\;k\in\N.
\]
And since $(x_k)$ was selected arbitrarily, such a sequence of
indices $\la_n$ exists for each choice of ``target-sequence''
$(x_k)\in V^\N$. This means that $\vthi'\in\Thi(A,\mf)$.

Thus, $\Thi(A,\mf)$ is dense in $A^\N$. \vspace*{0.2cm}

I will now show that $\Thi(A,\mf)$ is $G_\de$ in $A^\N$.

Indeed let $\big\{\X_\la\big\}^\infty_{\la=1}$ be a sequence of elements of $V^\N$ with finite support, that is dense in $V^\N$.

For each $\la=1,2,\ld$ I define
\[
\Thi_{\X_\la}=\left\{\begin{array}{l}
                       (\vthi)\in A^\N:\;\exists\;n_0\in\N,\;\text{such that}\;d([T^i(\vthi)]_{n_0},(\X_\la)_i)<\frac{1}{\la}  \\
                       \text{where $i$ runs from 1 to the last nonzero term of $\X_\la$}.
                     \end{array}\right\}
\]
Now, let $M$ be the index of the last nonzero  term of $\X_\la$. That means $(\X_\la)_i=0$ $\fa\;i>M$. Now,
\[
\Thi_{\X_\la}=\bigcap^M_{i=1}\bigg(\bigcup^\infty_{m=1}(T^i)^{-1}\bigg(V\times V\times\cdots\times V\times\underbrace{B\bigg((\X_\la)_i,\frac{1}{\la}\bigg)}_{m\text{-position}} \times V\times\cdots\bigg)\bigg).
\]
The  set $(T^i)^{-1}\Big(V\times V\times\cdots\times V\times
B\Big((\X_\la)_i,\frac{1}{\la}\Big)\times V\times\cdots\Big)$ is
open since it is the inverse image of an open set under a continuous
function $(T^i)$ in the product topology. Since $M<\infty$, then
$\Thi_{\X_\la}$ is a finite intersection of open sets, thus it is
open.

However, $\Thi(A,\mf)=\bigcap\limits^\infty_{\la=1}\Thi_{\X_\la}$; that means that $\Thi(A,\mf)$ is $G_\de$ in $A^\N$.

I have shown that $\Thi(A,\mf)$ is dense and $G_\de$ in $A^\N$.

Now, for the second part of Theorem \ref{thm4.1}, let $A$ be a vector space. I will deduce from the previous result that $\Thi(A,\mf)$ contains a linear subspace of $V^\N$, with the exception of 0, that is dense in $A^\N$.

Since $V$ is metrizable and separable, $V^\N$ is also metrizable and separable. Let $d^\ast$ be a metric in $V^\N$ that generates the product topology. For instance let $d^\ast$ be:
\[
d^\ast((a),(\vthi))=\sum^\infty_{n=1}\frac{1}{2^n}\frac{d(a_n,\vthi_n)}{1+d(a_n,\vthi_n)}.
\]
Let $\vthi'_k=(\vthi'_k)_n$ where $(\vthi'_k)_n\in A^\N$
$\fa\;n=1,2,\ld$  be a sequence of elements of $A^\N$ that is dense
in $A^\N$. Such a sequence exists, since $V^\N$ is separable and
metrizable.

Since $\Thi(A,\mf)$ is dense in $A^\N$, there exists $\vthi_1\in\Thi(A,\mf)$ such that $d^\ast(\vthi_1,\vthi'_1)<\frac{1}{1}$.

Since $\vthi_1\in\Thi(A,\mf)$, if I consider a sequence of targets
$(x_k)$: $x_k=0$ $k=1,2,\ld$ then there exists a sequence $\la_n$,
$(\la_n\in\mf$ $n=1,2,\ld)$ such that
$\dis\lim_{n\ra\infty}[T^k(\vthi_1)]_{\la_n}=0$ $k=1,2,\ld\;.$

I consider $\mf_1=\{\la_n\}$  to be that set of indices. $\mf_1$ is
infinite, thus it will be true that $\Thi(A,\mf_1)$ is dense in
$A^\N$.

Now  I inductively consider sequences $\vthi_{i+1}\in\Thi(A,\mf_i)$
such that $d^\ast(\vthi_{i+1},\vthi'_{i+1})<\frac{1}{i}$ and define
$\mf_{i+1}\subset\mf_i$ as the set of indices $\la'_n\in\mf_i$ such
that $\dis\lim_{n\ra\infty}[T^k(\vthi_{i+1})]_{\la'_n }=0$
$k=1,2,\ld\;.$

Clearly, the set of sequences $\{\vthi_i\}$ $i=1,2,\ld$ is dense in
$A^\N$. I will show that $\Thi(A,\mf)$ contains the linear subspace
of $V^\N$ generated by $\vthi_i$ $(i=1,2,\ld)$ with the exception of
0.

Let $a_1,\ld,a_\mi\in\R^\ast$. I need to show that $a_1\vthi_1+\cdots+a_\mi\vthi_\mi\in\Thi(A,\mf)$.

Let $(x_k)\in V^\N$ be a sequence of ``targets''. Since
$\vthi_\mi\in\Thi(A,\mf_{\mi-1})$, then there exist a sequence of
indices $\la_n\in\mf_{\mi-1}\subset\mf$ such that
$\dis\lim_{n\ra\infty}[T^k(\vthi_\mi)]_{\la_n}=x_k\cdot\frac{1}{a_\mi}$
$k=1,2,\ld\;.$

Since
$\la_n\in\mf_{\mi-1}\subset\mf_{\mi-2}\subset\cdots\subset\mf_1$ if
will be true that $\dis\lim_{n\ra\infty}[T^k(\vthi_i)]_{\la_n}=0$
$i=1,\ld,\mi-1$.

Thus
\[
\lim_{n\ra\infty}[T^k(a_1\vthi_1+\cdots+a_\mi\thi_\mi)]_{\la_n}=a_1\cdot0+\cdots+
a_\mi\cdot\frac{x_k}{a_\mi}=x_k,  \qquad k=1,2,\ld\;.
\]
Moreover, $a_1\vthi_1+\cdots+a_\mi\vthi_\mi\in A^\N$ since $A$ is a
linear subspace of $V$.

Thus $a_1\vthi_1+\cdots+a_\mi\vthi_\mi\in\Thi(A,\mf)$, which concludes the proof of Theorem \ref{thm4.1}.  \qb \vspace*{0.2cm}

Now, I proceed with the proof of Theorem \ref{thm4.2} \vspace*{0.2cm}\\
\noindent
{\bf Proof of Theorem \ref{thm4.2}.} I will prove Theorem
\ref{thm4.2} using the following two lemmas: First, however, we make
a definition to simplify notation:
\setcounter{Def}{2}
\begin{Def}\label{Def4.3}
Given the natural numbers $v,\la_1,\ld,\la_k$ we define
$\f_i\equiv\f_i(v,\la_1,\ld,\la_k)=\ssum^{\la_i}_{j=1}T^{k+1-i}_{(v+\la_1+\cdots+
\la_i,v+\cdots+\la_{i-1}+j)}$ $i=1,\ld,k$.
\end{Def}
\setcounter{obs}{3}
\begin{obs}\label{obs4.4}
For instance, if $\vthi_1,\ld,\vthi_v\in B\subset V$ and
$\vthi_{v+1},\ld,\vthi_{v+\cdots+\la_i}\in\Ga=\{\Ga_1,\ld,\Ga_\mi\}\subset
V$, a finite subset of $V$, then for each $i=1,2,\ld,k$
\[
[T^{k+1-i}(\vthi)]_{v+\cdots+\la_i}=L(i)+a_1\Ga_1+\cdots+a_\mi\Ga_\mi \ \ a_j\ge0 \ \ j=1,\ld,\mi
\]
where $\ssum^\mi_{j=1}a_j=\f_i$ and $L(i)\in\,$conv$(B\cup\{0\})$, by Definition \ref{Def4.3} and Proposition \ref{prop2.0}.

Specifically $a_1\Ga_1+\cdots+a_\mi\Ga_\mi$ is a convex combination of $\f_i\cdot\Ga_j$ $j=1,\ld,\mi$.
\end{obs}
\setcounter{lem}{4}
\begin{lem}\label{lem4.5}
Let $V,A$, $\e>0$ and $k\in\N$ be given as in Theorem \ref{thm4.2}, $v_1\in\N$ and $M^0\subseteq A$ a finite subset of $A$, also be given $(\sharp M^0<\infty)$.

Then there exist finite subsets $M^0\subseteq M^1\subseteq\cdots\subseteq M^k\subseteq A$ of $A$ $(\sharp M^k<\infty)$ and a natural number $m_0\in\N$ such that for each natural number $m\ge m_0$, there exists a partition of $m$: $v(m),\la_1(m),\ld,\la_k(m)\vdash m$ such that, if $\f_i(m)\equiv\f_i(v(m),\la_1(m),\ld,\la_k(m))$ as in Definition \ref{Def4.3}, then the following are true:
\begin{enumerate}
\item[(a)] conv$(\f_i(m)\cdot M^i)\oplus B\Big(0,\frac{\de\big(\frac{\e}{6}\big)}{4}\Big)\supseteq$conv$(M^{i-1}
    \oplus$conv$(-M^{i-1})$ $i=1,2,\ld,k$ \vspace*{-0.1cm}
\item[(b)] $\frac{\la_i(m)}{v(m)+\cdots+\la_i(m)}\cdot [5\|M^{i-1}\|_\rho+1]\le\frac{\e}{6_k}$ $\rho=1,2,\ld,N_{\e/3}$, $i=1,2,\ld,k$
\item[(c)] $\frac{2\cdot\|M^k\|_\rho}{v(m)}<\frac{\e}{6\cdot(\sharp M^k)}$, $v(m)>v_1$, $v(m)>\frac{m}{2}$ $\rho=1,\ld,N_{\e/3}$
\end{enumerate}
where if $M$ is a finite subset of $V$, then $\|M\|_\rho=\max\{\|x\|_\rho:\;x\in M\}$.
\end{lem}
\begin{Proof}
By Corollary \ref{cor2.6} and Definition \ref{Def4.3} it is known
that
\[
\f_i(m)\in\bigg[\bigg(\frac{\la_i(m)}{v(m)+\cdots+\la_i(m)}\bigg)^{k+1-i}\cdot
\frac{1}{(k+1-i)!},\frac{2\la_i(m)}{v(m)+\cdots+\la_i(m)}\bigg].
\]
By setting $\ga_i(m)=\frac{\la_i(m)}{v(m)+\cdots+\la_i(m)}$ (b) can be written as:\vspace*{0.2cm} \\
(b$'$) $\ga_i(m)\cdot[5\|M^{i-1}\|_\rho+1]\le\frac{\e}{6k}$ $\rho=1,\ld,N_{\e/3}$ $i=1,\ld,k$ \vspace*{0.2cm}\\
and the previous bound can be written as:
\setcounter{equation}{0}
\begin{eqnarray}
\f_i(m)\in\bigg[\ga_i(m)^{k+1-i}\cdot\frac{1}{(k+1-i)!}, \;2\ga_i(m)\bigg].  \label{sec4eq1}
\end{eqnarray}
I will show that we can specify intervals $[c_i,d_i]\subseteq[0,1]$
for $\ga_i(m)$, that are independent of $m$, such that, when
$\ga_i(m)$ belongs to $[c_i,d_i]$ $i=1,\ld,k$, then (a) and (b$'$)
hold. Then I will demonstrate how, given a suitable number $m\in\N$
and a $k$-tuple of intervals for $\ga_i(m)$, we can define a
partition of $m$: $v(m),\la_1(m),\ld,\la_k(m)\vdash m$ such that:
$\frac{\la_i(m)}{v(m)+\cdots+\la_i(m)}\in[c_i,d_i]$.

(a) can be written as:
\[
[\f_i(m)\cdot\text{conv}(M^i)]\oplus B\bigg(0,\frac{\de(\e/6)}{4}\bigg)\supseteq\text{conv}(M^{i-1})\oplus\text{conv}(-M^{i-1}).
\]

Moreover, $\f_i(m)\cdot B\Big(0,\frac{\de(\e/6)}{4}\Big)\subseteq
B\Big(0,\frac{\de(\e/6)}{4}\Big)$.

Indeed, $\f_i(m)<1$ by Definition \ref{Def4.3} and Proposition \ref{prop2.0} and if $d$ is the metric in $V$, $d(0,\la_x)\le d(0,x)$ when $0\le\la\le1$, since $V$ is locally convex.

Thus, it is enough to show that
\[
\f_i(m)\cdot\bigg[\text{conv}(M^i)\oplus B\bigg(0,\frac{\de(\e/6)}{4}\bigg)\bigg]\supseteq\text{conv}(M^{i-1})\oplus
\text{conv}(-M^{i-1}) \ \ \text{or}
\]
\[
\text{conv}(M^i)\oplus B\bigg(0,\frac{\de(\e/6)}{4}\bigg)\supseteq\frac{1}{\f_i(m)}\cdot
[\text{conv}(M^{i-1})\oplus\text{conv}(-M^{i-1})].
\]
However,
\[
\text{conv}(M^{i-1})\oplus\text{conv}(-M^{i-1})\subseteq(2M^{i-1}\cup(-2)M^{i-1})
\]
since, if $x\in\text{conv}(M^{i-1})\oplus\text{conv}(-M^{i-1})$ then there exist $\mi\in\N$, $\la_1,\ld,\la_\mi$, $\la'_1,\ld,\la'_\mi\in\R_+$ such that $\ssum^\mi_{i=1}\la_i=1=\ssum^\mi_{i=1}\la'_i$ and $X=\ssum^{\mi'}_{i=1}\la_iA_i-\ssum^\mi_{i=1}\la'_i\cdot A_i$ for some $A_1,\ld,A_\mi\in M^{i-1}$.

However,
\[
\sum^\mi_{i=1}\la_iA_i-\sum^\mi_{i=1}\la'_iA_i=\sum^\mi_{i=1}\frac{\la_i}{2}\cdot2A_i
+\sum^\mi_{i=1}\frac{\la_i}{2}(-2A_i).
\]
Thus $x\in\text{conv}(2M^{i-1}\cup(-2)M^{i-1})$ since,
$\ssum^\mi_{i=1}\frac{\la_i}{2}+\ssum^\mi_{i=1}\frac{\la'_i}{2}=1$.

Therefore, it is enough for (a) to show that
\[
\text{conv}(M^i)\oplus B\bigg(0,\frac{\de(\e/6)}{4}\bigg)\supseteq\frac{1}{\f_i(m)}\cdot\text{conv}
(2M^{i-1}\cup(-2)M^{i-1})
\]
and since $\f_i(m)\ge[\ga_i(m)]^{k+1-i}\cdot\frac{1}{(k+1-i)!}$ by (\ref{sec4eq1}) it is enough to prove: \vspace*{0.2cm} \\
(a$'$): $\text{conv}(M^i)\oplus B\Big(0,\frac{\de(\e/6)}{4}\Big)\supseteq\frac{1}{[\ga_i(m)]^{k+1-i}}\cdot(k+1-i)!\cdot
\text{conv}(2M^{i-1}\cup(-2)M^{i-1})$ $i=1,\ld,k$.

For a finite subset $M$ of $V$, let
$\|M\|_{\max(\e)}=\dis\max_\rho\{\|M\|_\rho:\;\rho=1,\ld,N_\e\}$.

If for each $i=1,2,\ld,k$ there is a set $M^i$ such that
\begin{align}
\text{conv}(M^i)\oplus B\bigg(0,\frac{\de(\e/6)}{4}\bigg)\supseteq&\bigg[\frac{k}{\e}(10\cdot\|M^{i-1}\|_{\max(\e/3)}+2)\bigg]^{k+1-i}\nonumber\\
&\cdot(k+1-i)!\cdot\text{conv}(2M^{i-1}\cup(-2)M^{i-1}). \label{sec4eq2}
\end{align}
Then if
\begin{eqnarray}
\ga_i(m)\in\bigg[\frac{\e}{k}\frac{1}{10\|M^{i-1}\|_{\max(\e/3)}+2},
\frac{\e}{k}\cdot\frac{1}{5\|M^{i-1}\|_{\max(\e/3)}+1}\bigg]  \label{sec4eq3}
\end{eqnarray}
then both (a$'$) and (b$'$) are satisfied. (The lower bound gives (a$'$) and the upper bound (b$'$)). Such sets $M^i$ $(i=1,\ld,k)$ as in (\ref{sec4eq2}) can easily be constructed. (By induction, assuming $M^{i-1}$ is known).
\setcounter{obs}{5}
\begin{obs}\label{obs4.6}
If $\De$ is dense in $V$ and $B_j\in V$ $j=1,2,\ld$ is a countable collection of vectors of $V$, then for every $\e>0$, there exist $B_j^\ast(\e)\in\De$ such that:
\[
\text{conv}(\{B^\ast_j(\e)\})\oplus B(0,\e)\supseteq\text{conv}(\{B_j\}).
\]
\end{obs}

Indeed, I choose $B^\ast_j(\e)\in\De$ such that $d(B^\ast_j(\e),B_j)<\frac{\e}{2^i}$. Now, if $\bi=\ssum^\mi_{i=1}\la_iB_i\in\text{conv}(\{B_j:j=1,2,\ld\})$ for some $\mi\in\N^\ast$, $\la_i\in[0,1]$, $\ssum^\mi_{i=1}\la_i=1)$, then I consider $\bi^\ast(\e)=\ssum^\mi_{i=1}\la_iB^\ast_i(\e)$.

Now
\[
d(\bi,\bi^\ast(\e))\!\le\!\sum^\mi_{i=1}d(0,\la_i(B_i-B^\ast_i(\e))\!\le\!\sum^\mi_{i=1}
d(0,B_i-B^\ast(\e))\!=\!\sum^\mi_{i=1}d(B_i,B_i^\ast(\e))\!<\!\sum^\mi_{i=1}\frac{\e}{2^i}=\e.
\]
Now, in order to construct sets $M^i$ $i=1,\ld,k$ that satisfy
(\ref{sec4eq2}), we see that if $M^{i-1}$ is known, since conv$A$ is
dense in $V$, Observation \ref{obs4.6} implies that there exist
$A^\ast_1,\ld,A^\ast_q\in\text{conv}A$ such that $q\le20\sharp
M^{i-1}$ and
\begin{align*}
&\text{conv}(\{A^\ast_i:\;i=1,\ld,q\})\oplus B\bigg(0,\frac{\de(\e/6)}{4}\bigg)\supseteq
\text{conv}\bigg(\bigg[\frac{k}{\e}(10\cdot\|M^{i-1}\|_{\max(\e/3)}+2)\bigg]^{k+1-i}\\
&
\cdot(k+1-i)!\cdot(2M^{i-1}\cup(-2)M^{i-1})\bigg)
\end{align*}

Observation \ref{obs4.6} can be used since the set whose convex hall we consider in the right part, is finite with at most $2\cdot(\sharp M^{i-1})$ elements.

This is exactly what is needed for the set $M^i$ in (\ref{sec4eq2}).

Now, since $A^\ast_1,\ld,A^\ast_q\in\text{conv}A$, there exist a number $\mi\in\N$ and $A_1,\ld,A_\mi\in A$ such that $A^\ast_j\in\text{conv}(\{A_1,\ld,A_\mi\})$ $j=1,\ld,q$. I define $M^i=\{A_1,\ld,A_\mi\}\cup M^{i-1}$.

$M^i$ is finite and $M^{i-1}\subseteq M^i$ and for each $\ga_i(m)$ as in (\ref{sec4eq3}), (a$'$) and (b$'$) are satisfied.

In this way, since the first set $M^0$ is given I have constructed sets $M^0\subseteq M^1\subseteq\cdots\subseteq M^k$, $\sharp M^k<\infty$ that satisfy (\ref{sec4eq2}) and intervals
\[
[c_i,d_i]=\bigg[\frac{\e}{k}\cdot\frac{1}{10\|M^{i-1}\|_{\max(\e/3)}+2},
\frac{\e}{k}\frac{1}{5\|M^{i-1}\|_{\max(\e/3)}+1}\bigg]
\]
such that for the sets $M^0,\ld,M^k$ and for arbitrary real numbers $\ga_i(m)$ such that $\ga_i(m)\in[c_i,d_i]$ (a$'$) and (b$'$) are true.

By the definition of the intervals $[c_i,d_i]$, we see that $d_1-c_1\ge d_2-c_2\ge\cdots\ge d_k-c_k$. Moreover $d_1+\cdots+d_k<\e$.

Let $m_0$ be a natural number such that $m_0>\frac{2}{d_k-c_k}$ and
$m_0>2\cdot v_1$ and $m_0>\frac{4}{3}\cdot 6(\sharp
M^k)\cdot\|M^k\|_\rho$ $\rho=1,\ld,N_{\e/3}$.

Now, for each $m\ge m_0$ I will construct a partition of $m$: $v(m),\la_1(m),\ld,\la_k(m)\vdash m$ such that $\frac{\la_i(m)}{v(m)+\cdots+\la_i(m)}\in[c_i,d_i]$.

Since $\frac{1}{m}<\frac{d_k-c_k}{2}$ there is a natural number $\la_k(m)\in\N$ such that $\frac{\la_k(m)}{m}\in[c_k,d_k]$.

Moreover, $m-\la_k(m)>\frac{m}{2}$ since
$\frac{m}{2}>\la_k(m)\Leftrightarrow\frac{1}{2}>\frac{\la_k(m)}{m}$
which is true since
$\frac{1}{2}>\e>d_1+\cdots+d_k>d_k>\frac{\la_k(m)}{m}$.

Thus $\frac{1}{m-\la_k(m)}<d_k-c_k<d_{k-1}-c_{k-1}$.

For each $i=k-1,k-2,\ld,1$, I consider the following proposition:

There are $\la_{i+1}(m),\ld,\la_k(m)$ such that
\begin{align*}
&\frac{\la_k(m)}{m}\in[c_k,d_k],\frac{\la_{k-1}(m)}{m-\la_k(m)}\in[c_{k-1},d_{k-1}],
\ld,\\
&\frac{\la_{i+1}(m)}{m-\la_k(m)-\cdots-\la_{i+2}(m)}\in[c_{i+1},d_{i+1}]
\end{align*}
and $m-\la_k(m)-\cdots-\la_{i+1}(m)>\frac{m}{2}$.

This proposition is true when $i=k-1$ as we showed earlier. I consider it to be true for $i\in\{2,\ld,k-1\}$ and will prove it for $i-1$: Since
\[
m-\la_k(m)-\cdots-\la_{i+1}(m)>\frac{m}{2}\Leftrightarrow\frac{1}{m-\la_k(m)-\cdots-
\la_{i+1}(m)}<\frac{2}{m}<d_k-c_k<d_i-c_i.
\]
Thus, there is a natural number $\la_i(m)$ such that $\frac{\la_i(m)}{m-\cdots-\la_{i+1}(m)}\in[c_i,d_i]$.

Moreover
\[
\frac{\la_k(m)+\cdots+\la_i(m)}{m}<d_i+\cdots+d_k<\e<\frac{1}{2}
\]
by the induction argument.

Thus,
\[
\frac{m}{2}>\la_k(m)+\cdots+\la_i(m)\Leftrightarrow m-\la_k(m)-\cdots-\la_i(m)>\frac{m}{2}.
\]
Thus, the finite induction is completed and I have constructed numbers $\la_1(m),\ld,\la_k(m)$ such that $\frac{\la_i(m)}{m-\la_k(m)-\cdots-\la_{i+1}(m)}\in[c_i,d_i]$ and $\la_1(m)+\cdots+\la_k(m)<\frac{m}{2}$.

I set $v(m)\!=\!m-\la_k(m)-\cdots-\la_1(m)$, and for the partition\\
$v(m),\la_1(m),\ld,\la_k(m)\!\vdash\! m$ of $m$ it is true that
\[
\frac{\la_i(m)}{v(m)+\la_1(m)+\cdots+\la_i(m)}=\frac{\la_i(m)}{m-\la_k(m)-\cdots-
\la_{i+1}(m)}\in[c_i,d_i].
\]
Moreover, by definition $v(m)>\frac{m}{2}$.

Thus, $v(m)>\frac{m_0}{2}$.

Thus, $v(m)>v_1$ and
$\frac{2\cdot\|M^k\|_\rho}{v(m)}<\frac{\e}{6(\sharp M^k)}$
$\rho=1,\ld,N_{\e/3}$.

That means that $v(m)$ satisfies (c).

Moreover,  since the sets $M^0,M^1,\ld,M^k$ and the real numbers
$\ga_i(m)=\frac{\la_i(m)}{v(m)+\cdots+\la_i(m)}$ satisfy (a$'$) and
(b$'$), then the sets $M^0,M^1,\ld,M^k$ and the natural numbers\\
$v(m),\la_1(m),\ld,\la_k(m)\vdash m$, $\fa\;m\ge m_0$ satisfy (a)
and (b).

This concludes the proof of Lemma \ref{lem4.5}.  \qb
\end{Proof}
\setcounter{lem}{6}
\begin{lem}\label{lem4.7}
Let $V,A,x_1,\ld,x_k$, $\e>0$ and $k\in\N$ be given as in Theorem
\ref{thm4.2}. Let $M^0\subseteq M^1\subseteq\cdots\subseteq M^k$ be
finite subsets of $A$ and $v,\la_1,\ld,\la_k$ natural numbers such
that if $\f_i=\f_i(v,\la_i,\ld,\la_n)$ as in Definition \ref{Def4.3}
and $m=v+\la_1+\cdots+\la_k$, then $M^0,\ld,M^k$,
$v,\la_1,\ld,\la_k,m$ satisfy (a) and (c) of Lemma \ref{lem4.5}.

If furthermore $M^0$ is such that $0,x_i\in\text{conv}(M^0)\oplus
B\Big(0,\frac{\de(\e/6)}{4}\Big)$ $i=1,\ld,k$ and if
$\vthi_1,\ld,\vthi_v\in V$ are given such that if $\vthi$ ``starts
with'' $\vthi_1,\ld,\vthi_v$ then,
$[T^i(\vthi)]_v\in\text{conv}(M^0)\oplus
B\Big(0,\frac{\de(\e/6)}{4k}\Big)$ $i=1,2,\ld,k$.

Then, there exist terms $\vthi_{v+1},\ld,\vthi_{v+\la_1+\cdots+\la_k}$ where:
\[
\begin{array}{ll}
  j=1,\ld,\la_1 & \vthi_{v+j}\in M^1 \\
  j=\la_1+1,\ld,\la_1+\la_2 & \vthi_{v+j}\in M^2 \\
  \multicolumn{2}{c}{\dotfill} \\
  \multicolumn{2}{c}{\dotfill}\\
  j=\la_1+\cdots+\la_{k-1}+1,\ld,\la_1+\cdots+\la_k & \vthi_{v+j}\in M^k
\end{array}
\]
such that if $\vthi\in V^\N$ ``starts with'' $\vthi_1,\ld,\vthi_{v+\la_1+\cdots+\la_k}$ then: \smallskip

(a)
$\|[T^{k+1-i}(\vthi)]_{v+\cdots+\la_i}-x_{k+1-i}\|_\rho<\frac{\e}{3}$
$\rho=1,\ld,N_{\e/3}$ $i=1,\ld,k$.

(b)
$\|[T^{k+1-i}(\vthi)]_{v+\cdots+\la_{i-1}+j}\|_\rho\le5\cdot\|M^{i-1}\|_\rho+1$
$\rho=1,\ld,N_{\e/3}$, $i=1,\ld,k$,\linebreak $j=1,\ld,\la_i$.
\end{lem}
\begin{Proof}
\setcounter{obs}{7}
\begin{obs}\label{obs4.8}
If $\thi_1,\ld,\thi_v$ are given as in Lemma \ref{lem4.7} and $\vthi_{v+1},\ld,\vthi_{v+a}\in M^i$ are $a$-many terms of $M^i$ $(a\in\N)$, then if $\vthi\in V^\N$ ``starts with'' $\vthi_1,\ld,\vthi_{v+a}$ it is true that:
\[
[T^j(\vthi)]_{v+a}\in\text{conv}(M^i)\oplus B\bigg(0,\frac{\de(\e/6)}{4}\bigg) \ \ j=1,\ld,k.
\]
\end{obs}

Indeed, by Corollary \ref{cor2.11} $[T^j(\vthi)]_{v+a}$ is a convex
combination of $[T^i(\vthi))]_v,\ld,$
$[T^j(\vthi)]_v,\vthi_{v+1},\ld,\vthi_{v+a}$. Now, since
$[T^i(\vthi)]_v\in\text{conv}(M^i)\oplus
B\Big(0,\frac{\de(\e/6)}{4}\Big)$ $i=1,\ld,k$ and $\vthi_{v+j}\in
M^i$ $j=1,\ld,a$, a convex combination of these elements will belong
to
\[
B\underbrace{\bigg(0,\frac{\de(\e/6)}{4k}\bigg)\oplus\cdots\oplus B\bigg(0,\frac{\de(\e/6)}{4k}\bigg)}_{k\;\text{times}}\oplus\text{conv}(M^i)\subseteq
B\bigg(0,\frac{\de(\e/6)}{4}\bigg)\oplus\text{conv}(M^i).
\]
Now, to prove Lemma \ref{lem4.7}, I will inductively define terms
$\vthi_{v+1},\ld,\vthi_{v+\la_1+\cdots+\la_k}$ in $k$-many steps: In
the $i^{\text{th}}$ step $(i\!=\!1,\ld,k)$ I define terms
$\vthi_{v+\cdots+\la_{i-1}+1},\ld,\vthi_{v+\cdots+\la_{i-1}+\la_i}\!\in\!
M^i$ as follows:

By Observation \ref{obs4.8} and the fact that in previous steps the terms belong to $M^{i-1}$, it is true that $[T^{k+1-i}(\vthi)]_{v+\cdots+\la_{i-1}}\in\text{conv}(M^{i-1})\oplus B\Big(0,\frac{\de(\e/6)}{4}\Big)$.

For the $1^{\text{st}}$ step this is given by the Lemma. (It is one of the conditions).

Temporarily, even if $0\notin A$, I set $\vthi_{v+\la_1+\cdots+\la_{i-1}+j}=0$ $j=1,\ld,\la_i$. Thus, for a sequence $\vthi\in V^\N$ that ``starts with'' $\vthi_1,\ld,\vthi_{v+\cdots+\la_i}$ I consider the vector $S=[T^{k+1-i}(\vthi)]_{v+\cdots+\la_i}$. We see that
\begin{eqnarray}
S=\sum^{v+\cdots+\la_{i-1}}_{j=1}T^{k+1-i}_{(v+\cdots+\la_i,j)}\cdot\vthi_j  \label{sec4eq4}
\end{eqnarray}
since  $\vthi_{v+\cdots+\la_{i-1}+j}=0$ for $j=1,\ld,\la_i$. Here
$T^i_{(n,m)}$ is the element of the $n^{\text{th}}$ row and
$m^{\text{th}}$ column of the matrix of $T^i$.

By Observation \ref{obs4.8}, $S$ is a convex combination of an element of conv$(M^{i-1})\oplus B\Big(0,\frac{\de(\e/6)}{4}\Big)$ and 0. However, $0\in\text{conv}(M^{i-1})\oplus B\Big(0,\frac{\de(\e/6)}{4}\Big)$ (given by the lemma).

Thus,
\begin{eqnarray}
S\in\text{conv}(M^{i-1})\oplus B\bigg(0,\frac{2\de(\e/6)}{4}\bigg).  \label{sec4eq5}
\end{eqnarray}
Now it is given that $x_{k+1-i}\in\text{conv}(M^0)\oplus
B\Big(0,\frac{\de(\e/6)}{4}\Big)\Rightarrow
x_{k+1-i}\in\text{conv}(M^{i-1})\oplus
B\Big(0,\frac{\de(\e/6)}{4}\Big)$. Thus,
\[
x_{k+1-i}-S\in\text{conv}(M^{i-1})\oplus\text{conv}(-M^{i-1})\oplus
B\bigg(0,\frac{3\de(\e/6)}{4}\bigg).
\]
Thus,
\[
\|x_{k+1-i}-S\|_\rho\le 2\|M^{i-1}\|_\rho+\frac{\e}{6}, \ \ \rho=1,\ld,N_{\e/3}.
\]
However, by Lemma \ref{lem4.5} it is true that
\[
\text{conv}(\f_iM^i)\oplus B\bigg(0,\frac{\de(\e/6)}{4}\bigg)\supseteq\text{conv}(M^{i-1})\oplus\text{conv}
(-M^{i-1}).
\]
Thus $x_{k+1-i}-S\in\text{conv}(\f_iM^i)\oplus
B\Big(0,\frac{4\de(\e/6)}{4}\Big)$.

This means that there is $x'_{k+1-i}\in V$ such that
$d(x'_{k+1-i},x_{k+1-i})<\de(\e/6)$. (Thus,
$\|x'_{k+1-i}-x_{k+1-i}\|_\rho<\frac{\e}{6}$
$\rho=1,2,\ld,N_{\e/3}$) and $x'_{k+1-i}-S\in\text{conv}(\f_iM^i)$.

By combining the previous two inequalities of the seminorm $\|\,\|_\rho$, it will be true that
\begin{eqnarray}
\|x'_{k+1-i}-S\|_\rho<2\cdot\|M^{i-1}\|_\rho+\frac{\e}{3}\qquad \rho=1,2,\ld,N_{\e/3}.  \label{sec4eq6}
\end{eqnarray}
Now, since $x'_{k+1-i}-S\in\text{conv}(\f_iM^i)$, there exist
$g_1,\ld,g_\mi\in\R_+$ and $M^i_1,\ld,M^i_\mi\in M^i$ such that
$\mi<\sharp(M^i)$, $\ssum^\mi_{j=1}g_j=1$,
$\ssum^\mi_{j=1}g_j\f_iM_j=x'_{k+1-i}-S$.

We will show that we can define terms
$\vthi_{v+\cdots+\la_{i-1}+1},\ld,\vthi_{v+\cdots+\la_i}\in\{M^i_1,\ld,M^i_\mi\}$
such that
\[
\sum^j_{\ti=1}T^{k+1-i}_{(v+\cdots+\la_i,v+\cdots+\la_{i-1}+\ti)}\cdot\vthi_{(v+\cdots+
\la_{i-1}+\ti)}\!=\!L_j\cdot(x'_{k+1-i}-S)+\ups_j: \; j\!=\!1,\ld,\la_i, \; L_j\in(0,1),
\]
for suitable $L_j$, $\ups_j$ such that
$\|\ups_j\|_\rho\le\frac{1}{3}$ $\rho=1,\ld,N_{\e/3}$.

Indeed, for each $N\in\N$, I can define $\vthi_{v+\cdots+\la_{i-1}+1},\ld,\vthi_{v+\cdots+\la_i}\in\{M^i_1,\ld,M^i_\mi\}$ such that $\fa\;j=1,\ld,\la_i$
\begin{eqnarray}
\sum^j_{\ti=1}T^{k+1-i}_{(v+\cdots+\la_i,v+\cdots+\la_{i-1}+\ti)}\cdot
\vthi_{(v+\cdots+\la_{i-1}+\ti)}=\frac{L}{N}(x'_{k+1-i}-S)+\bi_1M^i_1+\cdots+\bi_\mi M^i_\mi \label{sec4eq7}
\end{eqnarray}
for some $L\in\{1,\ld,N\}$ and some $\bi_\rho$ $\rho=1,\ld,\mi$ such
that $\bi_\rho<\max\Big\{\frac{g_\rho\f_i}{N},\frac{2}{v}\Big\}$.

(Of course $L,\bi_1,\ld,\bi_\mi$ depend on $j=1,\ld,\la_i$) and moreover for $j=\la_i$:
\begin{eqnarray}
\sum^{\la_i}_{\ti=1}T^{k+1-i}_{(v+\cdots+\la_i,v+\cdots+\la_{i-1}+\ti)}\cdot
\vthi_{(v+\cdots+\la_{i-1}+\ti)}=x'_{k+1-i}-S+\bi'_1M'_1+\cdots+\bi'_\mi\cdot M'_\mi  \label{sec4eq8}
\end{eqnarray}
where $|\bi'_i|<\frac{2}{v}$.

I can do this the following way

At the beginning I set  $\vthi_{v+\cdots+\la_{i-1}+\ti}=0$,
$\thi=1,\ld,\la_i$ as I mentioned before.

After that in $\la_i$ steps, I assign vectors to $\vthi_{v+\cdots+\la_{i-1}+\ti}$ from the set $\{M_1,\ld,M^i_\mi\}$.

For each $j=1,2,\ld,\la_i$  the sum
$\ssum^j_{\ti=1}T^{k+1-i}_{(v+\cdots+\la_i,v+\cdots+\la_{i-1}+\ti)}\cdot
\vthi_{v+\cdots+\la_{i-1}+\ti}$ will be of the form
$\ga_1(j)M^i_1+\cdots+\ga_\mi(j)M^i_\mi$ where
$\ga_1(j)+\cdots+\ga_\mi(j)<\f_i$ and
$\ga_1(\la_i)+\cdots+\ga_\mi(\la_i)=\f_i$.

After each step  the coefficient of some $M^i_j$ will increase by a
value less than $\frac{2}{v}$ since
$T^{k+1-i}_{(v+\cdots+\la_i,v+\cdots+\la_{i-1}+\ti)}<\frac{2}{v}$ by
Corollary \ref{cor2.2}.

In this way,  after $j_0$ (for some $1\le j_0\le\la_i)$ steps, I can
make the coefficient $\ga_1(j_0)$ of $M^i_1$ such that
$\ga_1(j_0)\in\Big(\frac{g_1\cdot\f_i}{N},\frac{g_1\cdot\f_i}{N}+\frac{2}{v}\Big)$.
After that I can do the same for the rest $M^i_2,\ld,M^i_\mi$ so
that after $j_1$ steps it is true that
$\ga_\ti(j_1)\in\Big(\frac{g_\ti\f_i}{N},\frac{g_\ti\f_i}{N}+\frac{2}{v}\Big)$
$\ti=1,\ld,\mi$ by assigning first the vector $M^i_2$, then $M^i_3$
etc. After I finish one ``round'' I can start assigning the vector
$M^i_1$ again so that after $j_2$ steps it is true that
\[
\ga_1(j_2)\in\bigg(\frac{2g_1\f_i}{N},\frac{2g_1\f_i}{N}+\frac{2}{v}\bigg) \ \ \text{and}
\]
\[
\ga_\ti(j_2)\in\bigg(\frac{g_\ti\f_i}{N},\frac{g_\ti\f_i}{N}+\frac{2}{v}\bigg) \ \ \ti=2,\ld,\mi
\]
and at the end of the second ``round'', it will be true that
\[
\ga_\ti(j_3)\in\bigg(\frac{2\cdot g_\ti\cdot\f_i}{N},\frac{2\cdot g_\ti\f_i}{N}+\frac{2}{v}\bigg) \ \ \text{for some} \ \ j_3 \ \ \ti=1,\ld,\mi.
\]
I continue like  this and at the end of the $N^{\text{th}}$ round I
assign vectors such that
$\ga_\ti(j_4)\in\Big(g_\ti\cdot\f_i-\frac{2}{v},g_\ti\cdot\f_i\Big)$
$\ti=1,\ld,\mi$ for some $j_4$.

Finally since the sum of the coefficients will be $\f_i$ at the
$(\la_i)^{\text{th}}$ step, I can assign vectors to the last terms
$\vthi_{v+\cdots+\la_i},\vthi_{v+\cdots+\la_i-1},\ld$ such that\\
$\ga_\ti(\la_i)\in\Big(g_\ti\cdot\f_i-\frac{2}{v},g_\ti\cdot\f_i+\frac{2}{v}\Big)$
$\ti=1,\ld,\mi$.

Now I have assigned vectors $M^i_1,\ld,M^i_\mi$ to the terms
$\vthi_{v+\cdots+\la_{i-1}+1},\ld,\vthi_{v+\cdots+\la_i}$ in such a
way that (\ref{sec4eq7}) and (\ref{sec4eq8}) hold.

From (\ref{sec4eq7}) and (\ref{sec4eq8}), I will now deduce (a) and (b) of the lemma.

(\ref{sec4eq4}) and (\ref{sec4eq8}) imply that
\[
[T^{k+1-i}(\vthi)]_{v+\cdots+\la_i}-x'_{k+1-i}=\bi'_1M^i_1+\cdots+\bi'_\mi M^i_\mi \ \ \text{where} \ \ |\bi'_1|<\frac{2}{v}.
\]
Thus
$\|[T^{k+1-i}(\vthi)]_{v+\cdots+\la_i}-x'_{k+1-i}\|_\rho<\mi\cdot\frac{2}{v}\cdot
\|M^i\|_\rho<\frac{2\cdot\|M^i\|_\rho\cdot(\sharp
M^i)}{v}<\frac{\e}{6}$ by Lemma \ref{lem4.5}: (c)
$\rho=1,\ld,N_{\e/3}$.

Moreover, we have shown that
$\|x_{k+1-i}-x'_{k+1-i}\|_\rho<\frac{\e}{6}$ $\rho=1,\ld,N_{\e/3}$.

Thus
\[
\|[T^{k+1-i}(\vthi)]_{v+\cdots+\la_i}-x_{k+1-i}\|_\rho<\frac{\e}{3} \ \ \rho=1,\ld,N_{\e/3} \ \ \text{which is (a)}.
\]
Now (\ref{sec4eq7}) implies that for each $j=1,\ld,\la_i$ and
$\rho=1,2,\ld,N_{\e/3}$
\[
\bigg\|\sum^j_{\ti=1}T^{k+1-i}_{(v+\cdots+\la_i,v+\cdots+\la_{i-1}+\ti)}\vthi_{(v+\cdots+
\la_{i-1}+\ti)}\bigg\|_\rho\le\frac{L}{N}\|x'_{k+1-i}-S\|_\rho+\bigg(\frac{\f_i}{N}+
\frac{2\mi}{V}\bigg)\cdot\|M^i\|_\rho
\]
where $L\in\{1,\ld,N\}$ and $N$ is the one I choose when I assign
vectors $M^i_1,\ld,M^i_\mi$ to $\vthi_{v+\cdots+\la_{i-1}+\ti}$
$\ti=1,\ld,\la_i$.

I can choose $N$ such that
$\Big(\frac{\f_i}{N}+\frac{2\mi}{V}\Big)\cdot\|M^i\|_\rho<\frac{1}{3}$
$\rho=1,2,\ld,N_{\e/3}$.

Indeed $\frac{2(\sharp M^i)\cdot\|M^i\|_\rho}{V}<\frac{\e}{6}$
$\rho=1,2,\ld,N_{\e/3}$ by Lemma \ref{lem4.5}: (c).

Thus, since by (\ref{sec4eq6})
$\|x'_{k+1-i}-S\|_\rho<2\|M^{i-1}\|_\rho+\frac{\e}{3}$,
$\rho=1,2,\ld,N_{\e/3}$ then
\begin{eqnarray}
\bigg\|\sum^j_{\ti=1}T^{k+1-i}_{(v+\cdots+\la_i,v+\cdots+\la_{i-1}+\ti)}\cdot
\vthi_{(v+\cdots+\la_{i-1}+\ti)}\bigg\|_\rho
\le2\|M^{i-1}\|_\rho+\frac{\e}{3}+\frac{1}{3}  \ \ \rho\!=\!1,\ld,N_{\e/3}.  \label{sec4eq9}
\end{eqnarray}
I consider now $[T^{k+1-i}(\vthi)]_{v+\cdots+\la_{i-1}+j}$ $j=1,\ld,\la_i$ and observe that
\[
[T^{k+1-i}(\vthi)]_{v+\cdots+\la_{i-1}+j}=S_{(j)}+\sum^j_{\ti=1}T^{k+1-i}_{(v+
\cdots+\la_{i-1}+j,v+\cdots+\la_{i-1}+\ti)}\cdot\vthi_{(v+\cdots+\la_{i-1}+\ti)}.
\]
where  for each $j=1,\ld,\la_i$, $S_{(j)}$ is a convex combination
of
$[T^1(\vthi)]_v,\ld,[T^{k+1-i}(\vthi)]_v,\linebreak\vthi_{v+1},\ld,\vthi_{v+\cdots+\la_{i-1}}$
and 0, by Corollary \ref{cor2.11}.

Thus, $S_{(j)}\in\text{conv}(M^{i-1})\oplus B\Big(0,\frac{2\de(\e/6)}{4}\Big)$ as in (\ref{sec4eq5}). Thus,
\begin{eqnarray}
\|S_{(j)}\|_\rho<\|M^{i-1}\|_\rho+\e/6 \qquad \rho=1,2,\ld,N_{\e/3}.  \label{sec4eq10}
\end{eqnarray}
Now, we observe that for each $j=1,\ld,\la_i$
\begin{align*}
[T^{k+1-i}(\vthi)]_{v+\cdots+\la_{i-1}+j}-S_{(j)}=&\sum^j_{\ti=1}
\frac{T^{k+1-i}_{(v+\cdots+\la_{i-1}+j,v+\cdots+\la_{i-1}+\ti)}}
{T^{k+1-i}_{(v+\cdots+\la_i,v+\cdots+\la_{i-1}+\ti)}} \\
&\cdot T^{k+1-i}_{(v+\cdots+\la_i,v+\cdots+\la_{i-1}+\ti)}\cdot\vthi_{(v+\cdots+
\la_{i-1}+\ti)}.
\end{align*}
By setting
\[
\al_\ti=\frac{T^{k+1-i}_{(v+\cdots+\la_{i-1}+j,v+\cdots+\la_{i-1}+\ti)}}
{T^{k+1-i}_{(v+\cdots+\la_i,v+\cdots+\la_{i-1}+\ti)}}
\]
and
\[
\bi_\ti=T^{k+1-i}_{(v+\cdots+\la_i,v+\cdots+\la_{i-1}+\ti)}\cdot\vthi_{(v+\cdots+\la_{i-1}+\ti)} \ \ \ti=1,\ld,j.
\]
I have that $[T^{k+1-i}(\vthi)]_{v+\cdots+\la_{i-1}+j}-S_{(j)}$ is
of the form $\ssum^j_{\ti=1}\al_\ti\bi_\ti$ where
$2\ge\al_1\ge\al_2\ge\cdots\ge\al_j>0$ by Proposition \ref{prop2.8}
and
$\Big\|\ssum^j_{\ti=1}\bi_\ti\Big\|_\rho<2\|M^{i-1}\|_\rho+\frac{\e}{3}+\frac{1}{3}$
$\rho=1,\ld,N_{\e/3}$ by (\ref{sec4eq9}).

Thus, by Proposition \ref{prop2.9}
$\Big\|\ssum^j_{\ti=1}\al_\ti\bi_\ti\Big\|\le4\cdot\|M^{i-1}\|_\rho+\frac{2}{3}+\frac{2\cdot\e}{3}$
$\rho=1,2,\ld,N_{\e/3}$.

This means, that for each $j=1,2,\ld,\la_i$
\[
\|[T^{k+1-i}(\vthi)]_{v+\cdots+\la_{i-1}+j}-S_{(j)}\|_\rho\le4\cdot\|M^{i-1}\|_\rho+\frac{2}{3}
+\frac{2\cdot\e}{3} \qquad \rho=1,2,\ld,N_{\e/3}.
\]
Thus, (\ref{sec4eq10}) implies now that for each $j=1,\ld,\la_i$
\[
\|[T^{k+1-i}(\vthi)]_{v+\cdots+\la_{i-1}+j}\|_\rho\le5\|M^{i-1}\|_\rho+\frac{2}{3}+\e
\le5\|M^{i-1}\|_\rho+1 \qquad \rho=1,2,\ld,N_{\e/3}
\]
which proves (b) and concludes the proof of Lemma \ref{lem4.7}.  \qb
\end{Proof}

I continue now with the proof of Theorem \ref{thm4.2}.

Let $a\in A$ an arbitrary element of $A$.

Since  $V$ is a locally convex metric space, the Cesaro operator $T$
preserves the limit of convergent sequences. Thus, if I consider the
sequence $\vthi'_n$: $\vthi'_i=\vthi_i$ $i=1,\ld,\rho$ and
$\vthi'_i=a$, $i\ge\rho+1$ where $\vthi_1,\ld,\vthi_\rho$ are the
terms given by Theorem \ref{thm4.2} then,
$\dis\lim_{v\ra\infty}[T^k(\vthi')]_v=a$ $k=1,2,\ld\;.$

Thus, there exists $v_1\in\N$ such that: If $\vthi_{\rho+1}=\cdots=\vthi_{v_1}=a$ and a sequence $\vthi'\in V^\N$ ``starts with'' $\vthi_1,\ld,\vthi_{v_1}$, then
\[
[T^i(\vthi')]_{v_1}\in B\bigg(a,\frac{\de(\e/6)}{4k}\bigg) \ \ i=1,\ld,k.
\]
\begin{obs}\label{obs4.9}
There exists a finite subset of $A$, $(M^0\subseteq A,\;\sharp
M^0<\infty)$ such that $0,x_1,\ld,x_k\in\text{conv}(M^0)\oplus
B\Big(0,\frac{\de(\e/6)}{4}\Big)$ and thus,
$\|x_i\|_\rho\le\|M^0\|_\rho+\frac{\e}{3}$: $i=1,2,\ld$,
$\rho=1,2,\ld,N_{\e/3}$ and $a\in M^0$.

Indeed, since $\overline{\text{conv}(A)}=V$, there exist
$0^\ast,x^\ast_1,\ld,x^\ast_k\in\text{conv}A$ such that
$d(x^\ast_i,x_i)<\frac{\de(\e/6)}{4}$ and
$d(0^\ast,0)<\frac{\de(\e/6)}{4}$.

Since $0^\ast,\ld,x^\ast_k\in\text{conv}A$ they can be expressed as
convex combinations of a finite number of elements of $A$. This
finite collection is the set $M^0$.

I mention that we can force $\|M^0\|_\rho\ge1$ $\rho=1,\ld,N_{\e/3}$
by adding a finite number of elements.

Now $V,A$, $\e>0$,  $k\in\N$, $M^0$ and $v_1$ satisfy the conditions
of Lemma \ref{lem4.5}, hence there exist sets $M^1,M^2,\ld,M^k$ and
a number $m_0$ as prescribed in Lemma \ref{lem4.5}.

We choose a number $m>m_0$ such that $m\in\mf$ and $m>2v_1$. By Lemma \ref{lem4.5} there is a partition $v,\la_1,\ld,\la_k\vdash m$ such that $v>v_1$ and conditions (a), (b), (c) of Lemma \ref{lem4.5} hold.

We also define terms $\vthi_{v_1+1},\ld,\vthi_v=a$. It will of course still be true that if a sequence $\vthi'\in V^\N$ ``starts with'' $\vthi_1,\ld,\vthi_v$, then
\[
[T^i(\vthi)]_v\in\text{conv}(M^0)\oplus B\bigg(0,\frac{\de(\e/6)}{4}\bigg) \ \ \text{since} \ \ a\in M^0
\ \ i=1,\ld,k.
\]

Now $V,A$, $\e>0$, $k\in\N$ $x_1,\ld,x_k$, $M^0,M^1,\ld,M^k$,
$v,\la_1,\ld,\la_k$ and $\vthi_1,\ld,\vthi_v$ all satisfy the
conditions of Lemma \ref{lem4.7}, thus there exist
$\vthi_{v+1},\ld,\vthi_{v+\la_1+\cdots+\la_k}\in M^k$ such that the
following hold: If $\vthi\in V^\N$ ``starts with''
$\vthi_1,\ld,\vthi_{v+\cdots+\la_k}$, then:
\begin{eqnarray}
\|[T^{k+1-i}(\vthi)]_{v+\cdots+\la_i}-x_{k+1-i}\|_\rho<\e/3 \ \ \rho=1,\ld,N_{\e/3} \ \ i=1,\ld,k  \label{sec4eq11}
\end{eqnarray}
and
\begin{align}
&\|[T^{k+1-i}(\vthi)]_{v+\cdots+\la_{i-1}+j}\|_\rho\le5\|M^{i-1}\|_\rho+1 \nonumber\\
&\rho=1,\ld,N_{\e/3}, \ \ i=1,\ld,k, \ \ j=1,\ld,\la_i.  \label{sec4eq12}
\end{align}
For $i=k$, (\ref{sec4eq11}) implies that
$\|[T(\vthi)]_{v+\cdots+\la_k}-x_1\|_\rho<\frac{\e}{3}$
$\rho=1,\ld,N_{\e/3}$. That means that
\begin{eqnarray}
d([T(\vthi)]_{v+\cdots+\la_k},x_1)<\e/3+\frac{\e}{3}<\e.  \label{sec4eq13}
\end{eqnarray}
For each $j=1,\ld,\la_k$ I consider
$[T^2(\vthi)]_{v+\cdots+\la_{k-1}+j}-x_2$. By Proposition
\ref{prop2.10} it is equal to
\[
=\frac{v\!+\!\cdots\!+\!\la_{k-1}}{v\!+\!\cdots\!+\!\la_{k-1}+j}[T^2(\vthi)]_{v\!+\!\cdots\!+\!\la_{k-1}}\!+\!
\frac{[T(\vthi)]_{v\!+\!\cdots\!+\!\la_{k-1}}\!+\!\cdots\!+\![T(\vthi)]_{v\!+\!\cdots\!+\!\la_{k-1}+j}}
{v+\cdots+\la_{k-1}+j}-x_2.
\]
Considering the $\|\,\|_\rho$ $\rho=1,2,\ld,N_{\e/3}$ seminorm
\begin{align*}
\|[T^2&(\vthi)]_{v+\cdots+\la_{k-1}+j}-x_2\|_\rho\le\frac{v+\cdots+\la_{k-1}}
{v+\cdots+\la_{k-1}+j}\cdot\|[T^2(\vthi)]_{v+\cdots+\la_{k-1}}-x_2\|_\rho\\
&+\frac{j}{v+\cdots+\la_{k-1}+j}\|x_2\|_\rho\\
&+\frac{\|[T(\vthi)]_{v+\cdots+\la_{k-1}+1}
\|_\rho+\cdots+\|[T(\vthi)]_{v+\cdots+\la_{k-1}+j}\|_\rho}{v+\cdots+\la_{k-1}+j}\\
\le\,&\e/3+\frac{j}{v+\cdots+\la_{k-1}+j}(\|x_2\|_\rho+5\cdot\|M^{k-1}\|_\rho+1),\;\text{by (\ref{sec4eq11}) and (\ref{sec4eq12})} \\
<\,&\e/3+\frac{\la_k}{v+\cdots+\la_k}(\|M^{k-1}\|_\rho+1+5\|M^{k-1}\|_\rho+1), \text{by Observation \ref{obs4.9}} \\
<\,&\e/3+2\cdot\frac{\e}{6k}<2\cdot\frac{\e}{3} \ \ \text{by Lemma \ref{lem4.5}: (b)}
\end{align*}
Thus for $i=1$ it is true that
\[
\|[T^{i+1}(\vthi)]_{v+\cdots+\la_{k-1}+j}-x_{i+1}\|_\rho\le\frac{2\e}{3} \ \ j=1,\ld,(\la_{k-i+1}+\cdots+\la_k) \ \ \rho=1,2,\ld,N_{\e/3}.
\]
For $i+1$: We consider
$[T^{i+2}(\vthi)]_{(v+\cdots+\la_{k-i-1}+j)}-x_{i+2}$:
$j=1,\ld,(\la_{k-i}+\cdots+\la_k)$.

By Proposition \ref{prop2.10} it is equal to:
\begin{align*}
&\frac{v+\cdots+\la_{k-i-1}}{v+\cdots+\la_{k-i-1}+j}\cdot[T^{i+2}(\vthi)]_{(v+\cdots+
\la_{k-i-1})}\\
&+\frac{[T^{i+1}(\vthi)]_{(v+\cdots+\la_{k-i-1}+1)}+\cdots+[T^{i+1}(\vthi)]_{(v+\cdots+
\la_{k-i-1}+j)}}{v+\cdots+\la_{k-i-1}+j}-x_{i+2}.
\end{align*}
Thus considering the $\|\,\|_\rho$ seminorm, $\rho=1,\ld,N_{\e/3}$
\begin{align*}
\|[T^{i+2}(\vthi)]_{(v+\cdots+\la_{k-i-1}+j)}\!-\!x_{i+2}\|_\rho\!\le\!&
\frac{v+\cdots+\la_{k-i-1}}{v+\cdots+\la_{k-i-1}+j}\!\cdot\!\|[T^{i+2}(\vthi)]_{v+\cdots+
\la_{k-i-1}}\!\!-\!\!x_{i+2}\|_\rho \\
&+\frac{j}{v+\cdots+\la_{k-i-1}+j}\cdot\|x_{i+2}\|_\rho \\
&+\frac{\la_{k-i}}{v+\cdots+\la_{k-i-1}+j}\cdot(5\cdot\|M^{k-i-1}\|_\rho+1)\\
&+
\frac{|j-\la_{k-i}|}{v+\cdots+\la_{k-i-1}+j}\bigg(\|x_{i+1}\|_\rho+\frac{2\e}{3}\bigg)
\end{align*}
by the inductive argument and (\ref{sec4eq12}) (where by considering
the absolute value of $(j-\la_{k-i})$ we cover the cases $1\le
j\le\la_{k-i}$ and $j>\la_{k-i}$ at the same time).
\begin{align*}
\le&\frac{v+\cdots+\la_{k-i-1}}{v+\cdots+\la_{k-i-1}+j}\cdot\frac{\e}{3}+\frac{k\cdot\e}{6k}
+\frac{\e}{6k}+
\frac{|j-\la_{k-i}|}{v+\cdots+\la_{k-i-1}+j}(\|M^{i-1}\|_\rho+1) \\
&\text{by (\ref{sec4eq11}) and Observation \ref{obs4.9} and Lemma \ref{lem4.5}: (b) }\\
\le&\frac{\e}{3}+\frac{\e}{6}+\frac{\e}{6k}+\frac{(k-1)\cdot\e}{6k}=\frac{2\e}{3} \ \ \text{by
Lemma \ref{lem4.5}: (b)}
\end{align*}
which completes the induction.

Thus, considering the index $(v+\cdots+\la_k)$ we will have
\[
\|[T^{i+1}(\vthi)]_{(v+\cdots+\la_k)}-x_{i+1}\|_\rho\le\frac{2\e}{3} \ \ \rho=1,\ld,N_{\e/3}.
\]
Thus
\[
d([T^{i+1}(\vthi)]_{(v+\cdots+\la_k)},x_{i+1})<\e \ \ i=1,\ld,k-1.
\]
Combining this with (\ref{sec4eq13}) we conclude the proof of Theorem \ref{thm4.2}.  \qb
\end{obs}
\setcounter{rem}{9}
\begin{rem}\label{rem4.10}
Theorem \ref{thm4.1} can be strengthened with a very slight modification of the proof as follows:

If $\Fi_n\in(R_+\cup\{\infty\})^\N$ such that $\dis\lim_{n\ra\infty}\Fi_n=\infty$ and $\Fi_n>\inf\{d(0,x):x\in A\}$ $\fa\;n$, then the conclusion of Theorem \ref{thm4.1} concerning the residuality holds with the modification that $\Thi(A,\mf)$ will be replaced by:
\[
\Thi(A,\mf,\{\Fi_n\})=\left\{\begin{array}{l}
                               (\vthi)\!\!\in\!\! A^\N\!:\!d(0,\vthi_n)\!<\!\Fi_n \ \ \fa\;n=1,2,\ld, \;\text{and such that}\:\fa(x_k)\!\in\! V^\N, \\
\;\text{sequence of ``targets'' there exists a sequence of indices} \; \la_n,\\
\la_n\in\mf,\;\text{such that}\;[T^k(\vthi)]_{\la_n}\!\!\ra\!\! x_k \;\text{as}\; \la_n\!\!\ra\!\!\infty\;k=1,2,\ld\;.\\
                             \end{array}\right.
\]
Furthermore, if $A$ is a vector space and $\Fi_n=+\infty$ $\fa_n$, then $\Thi(A,\mf)=\Thi(A,\mf,\{\Fi_n\})$ contains a dense vector subspace, with the exception of 0.
\end{rem}
\begin{rem}\label{rem4.11}
The results of \S\,4 imply easily that in the case when $A$ is a vector space, then the set $\Thi_{A,V}$ of Theorem \ref{thm3.1} contains a dense in $A^\N$ vector subspace, with the exception of 0.
\end{rem}

Finally, we show that if $A$ is bounded then $\Thi(A,\mf)=\emptyset$. This follows easily from Proposition \ref{prop4.12} below.
\setcounter{prop}{11}
\begin{prop}\label{prop4.12}
Let $A\subset[0,1]$ and $a\in A^\N$ such that for some $n\in\N$ we have $[T(a)]_n<\frac{1}{8}$. Then $[T^2(a)]_n<\frac{15}{16}$.
\end{prop}
\begin{Proof}
Since $[T(a)]_n<\frac{1}{8}$ then at least half of the terms
$a_1,\ld,a_n$ are: $a_i<\frac{1}{4}$. Now,\vspace*{-0.2cm}
\[
[T^2(a)]_n=\sum^\infty_{i=1}T^2_{(n,i)}\cdot a_i=\sum^n_{i=1}T^2_{(n,i)}\cdot a_i.\vspace*{-0.2cm}
\]
However, by Proposition \ref{prop2.1} $T^2_{(n,i)}\ge
T^2_{(n,i+1)}$.

Thus, $[T^2(a)]_n\le\sum^n_{i=1}T^2_{(n,i)}\cdot a'_i$ where
$a'_1,\ld,a'_n$ is a rearrangement of $a_1,\ld,a_n$ so that it is
decreasing $(a'_i\ge a'_{i+1})$.

We know that $a'_1,\ld,a'_{n/2}\le1$ and that
$a'_{\frac{n}{2}+1},\ld,a'_n<\frac{1}{4}$.

Thus,\vspace*{-0.2cm}
\[
[T^2(a)]_n\le\sum^{n/2}_{i=1}T^2_{(n,i)}+\frac{1}{4}\cdot\sum^n_{i=n/2}T^2_{(n,i)}.\vspace*{-0.2cm}
\]
But by Corollary \ref{cor2.6}
\[
\sum^n_{i=n/2}T^2_{(n,i)}\ge\bigg(\frac{\frac{n}{2}}{n}\bigg)^2\cdot\frac{1}{2}=\frac{1}{8}.
\]
Thus
\[
[T^2(a)]_n\le\frac{7}{8}\cdot1+\frac{1}{8}\cdot\frac{1}{4}<\frac{15}{16}. \ \ \text{\qb}
\]
\end{Proof}
\bigskip
\noindent
{\bf Acknowledgement:}

I would like to thank my Professor V. Nestoridis whose assistance, throughout writing this paper was without doubt a determinative factor; in suggesting new directions, but also helping me with old techniques he played a crucial role in the expansion of this paper.
 \vspace*{0.5cm}
 University of Athens \vspace*{-0.2cm}\\
 Department of Mathematics  \vspace*{-0.2cm}\\
 157 84 Panepistimioupolis \vspace*{-0.2cm}\\
 Athens, GREECE\medskip\\
Douvropoulos Theodosios: douvr001@math.umn.edu, atheodosios52@gmail.com

\end{document}